\documentclass[journal ]{new-aiaa}
\usepackage[utf8]{inputenc}
\usepackage{textcomp}

\def\rebuttal

\usepackage{graphicx}
\usepackage[dvipsnames]{xcolor}
\usepackage[normalem]{ulem}
\usepackage{mathtools}   
\usepackage{acronym}
\usepackage{tikz}
\usepackage{pgfplots}
\pgfplotsset{compat=1.18}
\usepackage{array}
\usepackage{booktabs}
\usepackage{subcaption}
\usepackage{hyperref}
\usepackage[capitalize]{cleveref}
\crefname{assumption}{Assumption}{assumptions}
\usepackage{placeins}
\usepackage{url}
\usepackage{float}
\usepackage{multirow}
\usepackage[symbols,nogroupskip,sort=standard,stylemods=mcols,style=mcolindex]{glossaries-extra}
\usepackage{chngcntr}
\usepackage{apptools}
\usetikzlibrary{positioning, arrows.meta, shapes.multipart} 
\usepackage{rotating} 
\usepackage[symbols,nogroupskip,sort=standard,stylemods=mcols,style=mcolindex]{glossaries-extra}

\usepackage[version=4]{mhchem}
\usepackage{siunitx}
\usepackage{longtable,tabularx}
\usepackage{amsmath, amsfonts, amssymb, amsthm}

\ifdefined\rebuttal
    \newcommand{\remove}[1]{\textcolor{red}{\sout{#1}}}
    
\else
    \newcommand{\remove}[1]{}
    
\fi

\def\p{\mathrm{p}}

\newcommand{\Om}{\Omega}
\newcommand{\om}{\omega}

\newacronym{amr}{AMR}{Area-to-Mass Ratio}
\newacronym{aop}{AoP}{Argument of Perigee}
\newacronym{cnes}{CNES}{Centre National d'Études Spatiales}
\newacronym{em}{EM}{Expectation-Maximization}
\newacronym{gmm}{GMM}{Gaussian Mixture Model}
\newacronym{hpc}{HPC}{High-Performance Computing}
\newacronym{iid}{i.i.d.}{identically and independently distributed}
\newacronym{mc}{MC}{Monte Carlo}
\newacronym{od}{OD}{Orbit Determination}
\newacronym{pca}{PCA}{Principal Component Analysis}
\newacronym{pdf}{p.d.f.}{probablity distribution function}
\newacronym{raan}{RAAN}{Right Ascension of the Ascending Node}
\newacronym{ssa}{SSA}{Space Situational Awareness}
\newacronym{wrt}{w.r.t.}{with respect to}

\DeclareMathAlphabet{\mathmybb}{U}{bbold}{m}{n}

\newcounter{inlineenum}
\renewcommand{\theinlineenum}{\alph{inlineenum}}
\newenvironment{inlineenum}
{\unskip\ignorespaces\setcounter{inlineenum}{0}%
\renewcommand{\item}{\refstePoCounter{inlineenum}{\textit{\theinlineenum})~}}}
{\ignorespacesafterend}

    \acrodef{cnes}[CNES]{Centre National d'Études Spatiales}
    
    \acrodef{cdm}[CDM]{Collision Data Message}
    
    \acrodef{geo}[GEO]{Geosynchronous Orbit}
    
    \acrodef{hbr}[HBR]{Hard-Body Radius}
    
    \acrodef{isae}[ISAE-SUPAERO]{Institut Supérieur de l'Aéronautique et de l'Espace}
    
    \acrodef{leo}[LEO]{Low-Earth Orbit}
    
    \acrodef{poc}[PoC]{Probability of Collision}
    
    \acrodef{resp}[resp.]{respectively}
    
    \acrodef{spoc}[sPoC]{Scaled Probability of Collision}
    
    \acrodef{tca}[TCA]{Time of Closest Approch}

\begin{document}

\def\ourtitle{A Physics-Informed Statistical Learning Model for Long-Term Fragmentation Cloud Propagation}

\def\NameYema{Yema \textsc{Paul}}

\def\DivSME{DOA/SME/SE}
\def\DivOR{DTN/DV/OR}
\def\OrgCNES{\acs{cnes}}
\def\StrCNES{18 Av. Édouard Belin}
\def\ZiPoCNES{31400}

\def\DivLab{Space Advanced Concepts Laboratory}
\def\OrgISAE{\acs{isae}}
\def\StrISAE{10 Av. Édouard Belin}
\def\ZipISAE{31400}

\def\CityToul{Toulouse}
\def\CountryFR{France}

\def\JobYema{PhD candidate, \DivLab}
\def\JobPau{Head of orbit determination office, \DivOR}

\author{
Yema Paul\footnote{Ph.D. Candidate, ISAE-SUPAERO, Toulouse, France. Corresponding author: yema.paul@isae-supaero.fr.},
Diogo Sousa\footnote{Master Student, ISAE-SUPAERO, Toulouse, France.},
Albert Servitje\footnote{Master Student, ISAE-SUPAERO, Toulouse, France.},
Rodrigo Pariente\footnote{Master Student, ISAE-SUPAERO, Toulouse, France.},
Emmanuel Delande\footnote{Space Surveillance Specialist, CNES, Toulouse, France.},
Christophe Taillan\footnote{Space Surveillance Specialist, CNES, Toulouse, France.},
Joan-Pau Sanchez\footnote{Professor of Astrodynamics and Mission Design, ISAE-SUPAERO, Toulouse, France.}
}
\affil{ISAE-SUPAERO, Toulouse, France; CNES, Toulouse, France}

\def\ourabstract{
}

\def\ourkeywords{
Orbital fragmentation;
Space debris;
Generative density modeling;
Surrogate modeling;
Machine learning;
Gaussian mixture models}

\title{\ourtitle}

\maketitle


\acresetall




\section{Introduction}
\label{sec:introduction}

The long-term evolution of orbital fragmentation clouds is a central problem in space situational awareness, debris-environment modelling, and collision-risk assessment. Following a breakup event, hundreds of thousands of fragments may be generated and subsequently evolve under the combined effects of atmospheric drag, Earth oblateness, luni-solar perturbations, and solar radiation pressure. Predicting the future distribution of these fragments is essential for evaluating collision risk, assessing the impact of fragmentation events on the orbital environment, and supporting long-term sustainability studies \cite{Johnson_NL_2001_1,Anselmo_L_2009_1,giudici2023probabilistic,giudici2024density,Johnson_NL_2008_1}.

High-fidelity Monte Carlo simulations remain the reference approach for fragmentation-cloud modelling. Starting from an initial population generated by statistical breakup models such as the NASA Standard Breakup Model \cite{Johnson_NL_2001_1}, each fragment is propagated individually using a detailed orbital propagator. While highly accurate, this approach becomes computationally expensive for long-term studies involving hundreds of thousands of fragments and multiple fragmentation scenarios\cite{giudici2023probabilistic,paul24fragmentation}.

To alleviate this burden, several continuum and density-based formulations have been proposed, in which the debris cloud is represented by a probability density rather than a collection of individual fragments \cite{McInnes1993,Letizia_F_2016_1,Frey_S_2021_1,giudici2023probabilistic,giudici2024density}. These approaches propagate the cloud density by solving a continuity equation in orbital-element space, typically through the Method of Characteristics.
To reduce the dimensionality of the problem, many models rely on the classical band-formation assumption. As differential orbital precession disperses the cloud, the angular variables are assumed to become uniformly distributed and therefore independent of the remaining orbital elements, reducing the density representation from six to three dimensions \cite{Letizia_F_2016_1,Frey_S_2021_1}.

This approximation yields substantial computational savings. However, recent studies suggest that band formation may require several decades to be achieved for highly inclined orbital regimes \cite{Paul2025FragmentationReport,Bazenga2024BandFormation,ParienteServitje2025CloudCollision}. Moreover, uniformity of a variable does not, in general, imply statistical independence from the remaining variables \cite{Nelsen2006Copulas}. In particular, numerical simulations show that the \gls{raan} may retain a strong dependence on the semi-major axis even after its marginal distribution becomes approximately uniform, as illustrated in Fig.~\ref{fig:raan_2y}. The present work therefore retains these dependencies through conditional density modeling rather than assuming independence of the angular variables.

In previous works \cite{paul24fragmentation,Taillan_C_2024_1b,Delande_E_2024_1}, the authors introduced point-process representations of fragmentation clouds that significantly reduced the number of propagated samples required for long-term population studies.

The present work formulates fragmentation-cloud propagation as a density-learning problem in Keplerian-element space and introduces a Hierarchical Generative Density Model (HGDM) that exploits the semi-major axis as the principal driver of long-term cloud evolution. Trained on a limited number of propagated representative fragments, the resulting surrogate provides a compact generative representation of the debris cloud capable of reconstructing high-fidelity multidimensional distributions while substantially reducing computational cost and storage requirements relative to conventional Monte Carlo propagation.\footnote{The methodology and validation results presented in this paper were previously reported at the International Conjunction Assessment Workshop \cite{Paul2025ICAWFragmentation} and in a CNES research note \cite{Paul2025FragmentationReport}. The present manuscript provides the corresponding archival journal version.}

The main contributions of this work are twofold: (i) the introduction of a physics-informed hierarchical generative density model for long-term orbital fragmentation-cloud propagation; and (ii) the demonstration that accurate cloud distributions can be reconstructed from a small set of propagated representative fragments, yielding order-of-magnitude reductions in computational cost and storage requirements.

\section{Problem Formulation}
\label{sec:problem_formulation}

Consider a fragmentation event occurring at epoch $t_0$. The objective of this work is to model the spatial distribution of the resulting debris cloud at future epochs $t\ge t_0$.
The cloud is represented in the Keplerian-element space
$
\mathcal K \subset \mathbb R^6,
$
whose coordinates are the semi-major axis $a$, eccentricity $e$, inclination $i$, \gls{raan} $\Omega$, \gls{aop} $\omega$, and mean anomaly $M$. These variables characterize the orbital geometry relevant for most debris-density and collision-risk applications.

Let
$
X_t=(a,e,i,\Omega,\omega,M)
$
denote the random state of a debris fragment at epoch $t$. The initial debris population is generated by a probabilistic breakup model in Cartesian space, namely the NASA Standard Breakup Model \cite{Johnson_NL_2001_1}. The resulting distribution can be mapped into Keplerian-element space \cite{frey2021transformation}, yielding an initial probability density function $p_0^{\mathrm{ref}}$ on $\mathcal K$.
Assuming that debris fragments evolve independently after the breakup and are subject to the same orbital dynamics, the debris cloud at any future epoch $t$ can be characterized by a probability density function $p_t^{\mathrm{ref}}$ on $\mathcal K$. In practice, $p_t^{\mathrm{ref}}$ is unknown and can only be approximated through the propagation of a sufficiently large Monte Carlo population. Throughout this paper, $p_t^{\mathrm{ref}}$ will be referred to as the reference debris-cloud distribution.
The objective of this work is to learn a surrogate density $p_k$ that approximates $p^{\mathrm{ref}}_k$ at each evaluation epoch. The learning problem is therefore not to predict the trajectory of a single fragment, but to reconstruct the probability law of the propagated ensemble
$p_k(x) \simeq p^{\mathrm{ref}}_k(x),
$ for $ x \in \mathcal K .$
In practice, the reference distribution $p_t^{\mathrm{ref}}$ is not available analytically and must be inferred from a finite set of propagated representative fragments. Let
$\mathcal D_k=
\left\{
x_k^{(j)}
\right\}_{j=1}^{N_{\mathrm{train}}}
$
denote the set of debris states obtained by propagating $N_{\mathrm{train}}$ representative fragments from the breakup epoch $t_0$ to the evaluation epoch $t_k$. The objective is to use $\mathcal D_k$ to construct a surrogate density $p_k$ that accurately approximates the reference distribution $p_k^{\mathrm{ref}}$.

Unlike conventional Monte Carlo approaches, which require storing and propagating a large population of debris fragments, the proposed framework seeks to learn a compact parametric representation of the cloud distribution. Once fitted, the surrogate model can efficiently generate synthetic debris-cloud realizations at epoch $t_k$ while preserving the dominant statistical properties of the reference population.
The problem considered in this work can therefore be viewed as a physics-informed generative density-estimation problem: high-fidelity orbital propagation provides representative samples of the debris cloud, while the surrogate model learns a compact probabilistic representation of the corresponding distribution.

\section{Hierarchical Generative Density Model}
\label{sec:HGDM}

The objective of the proposed framework is to construct a compact surrogate representation of the reference debris-cloud distribution $p_k^{\mathrm{ref}}$. A direct approximation of the full six-dimensional density would require a large number of parameters and a correspondingly large training dataset, making the resulting model difficult to estimate and sample from efficiently. Instead, we seek a structured representation that captures the dominant statistical dependencies while remaining computationally tractable.

The proposed approach is motivated by the observation that the semi-major axis largely governs the long-term evolution of fragmentation clouds. Differential orbital perturbations induce characteristic patterns linking the semi-major axis to other orbital elements, particularly the eccentricity, \gls{raan}, and \gls{aop}. As shown later in Section~\ref{sec:results}, these dependencies remain significant even after long propagation times and cannot be fully captured by the classical band-formation assumption (see \cref{fig:raan_2y}).
Accordingly, the debris-cloud distribution is represented through the hierarchical factorization
\begin{align}
p_k(a,e,i,\Omega,\omega,M) = p_k(a)\,p_k(i)\,p_k(e,\Omega,\omega|a)\,
p_k(M).
\label{eq:hierarchical_factorization}
\end{align}
This formulation can be interpreted as a hierarchical generative model in which the semi-major axis acts as the primary latent driver of the cloud morphology. Rather than estimating a generic six-dimensional density, the problem is reduced to estimating a marginal distribution, two independent distributions, and a lower-dimensional conditional density \cite{Paul2025ICAWFragmentation}. This substantially reduces model complexity while preserving the dominant structures observed in propagated debris clouds.
The following subsections describe the estimating procedure adopted for each component of Eq.~(\ref{eq:hierarchical_factorization}).

\subsection{Semi-Major-Axis Distribution}
\label{sec:sma}
The semi-major axis plays a central role in the formation and long-term evolution of fragmentation clouds. At fragmentation epoch, the velocity increments generated by the breakup induce instantaneous variations of the orbital elements through the impulsive form of Gauss' variational equations. Since the semi-major axis is directly related to orbital energy, it captures a large fraction of the dispersion introduced by the fragmentation process. Subsequently, the orbital elements evolve according to the continuous form of Gauss' variational equations under the action of orbital perturbations. Because many secular drift rates and perturbation effects depend strongly on the semi-major axis, significant statistical dependencies emerge between $a$ and other orbital elements, particularly the eccentricity, \gls{raan}, and \gls{aop}. For this reason, the semi-major axis is selected as the root variable of the hierarchical generative model introduced in Eq.~(\ref{eq:hierarchical_factorization}).

\begin{figure}[H]
	\centering
	\includegraphics[width=0.6\linewidth]
    {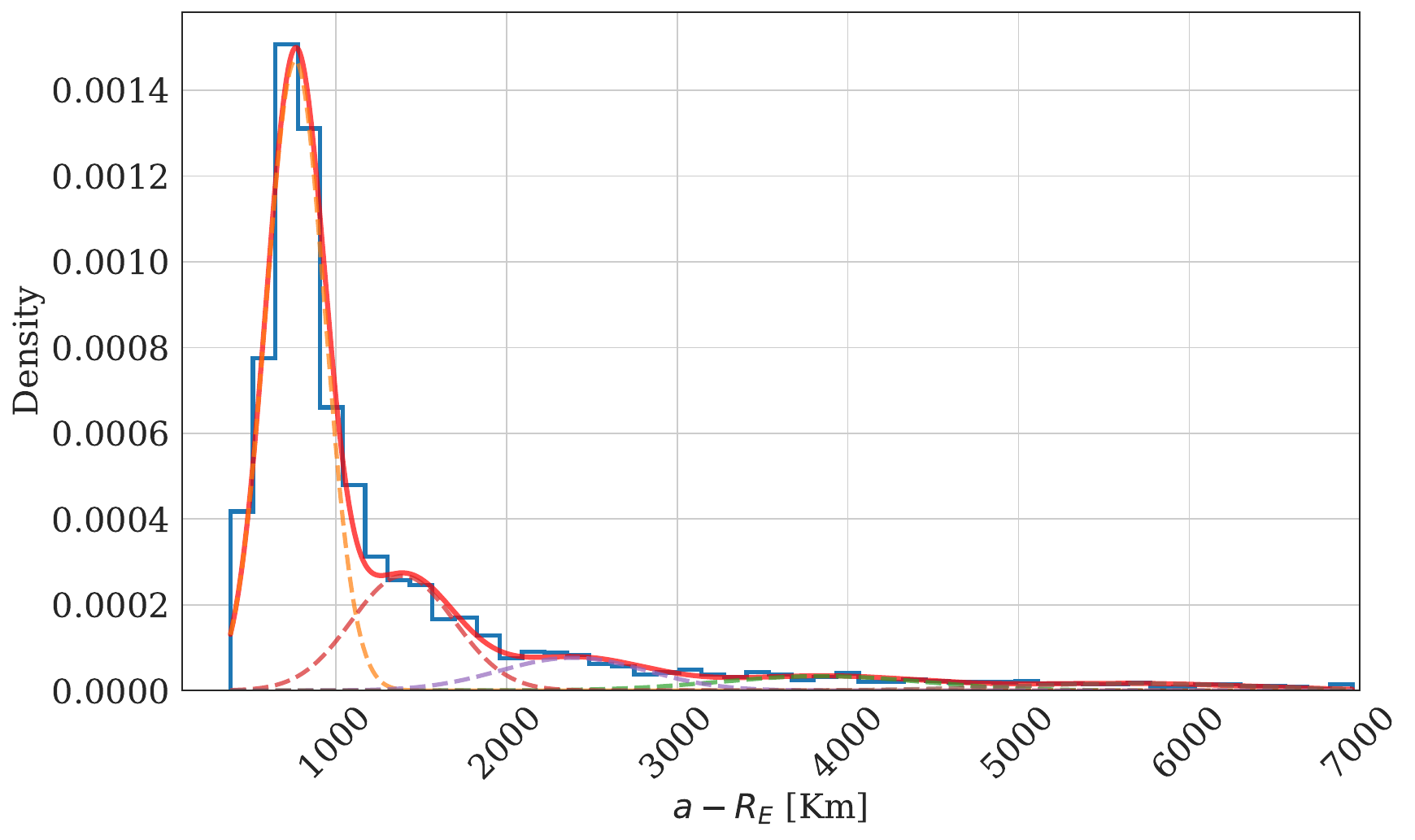}
\caption{GMM approximation of the semi-major-axis distribution two years after fragmentation epoch.}
    \label{fig:envisat_sma}
\end{figure}

The first step of the estimation procedure therefore consists in estimating the marginal density $p_k(a)$ from the propagated representative fragments. Inspection of the reference distributions reveals that the semi-major-axis density is generally multi-modal and may exhibit significant skewness and local concentration regions. A single parametric distribution is therefore insufficient to accurately represent the observed variability \cite{Letizia_F_2016_1}.
To accommodate these features, the semi-major-axis density is modeled using a Gaussian Mixture Model (GMM),
\begin{align}
p_k(a) = 
\sum_{m=1}^{N_G}
w_m\,
\mathcal N(a;\mu_m,\sigma_m^2),
\label{eq:gmm_a}
\end{align}
where $N_G$ is the number of mixture components, $w_m$ are non-negative weights satisfying
\(
\sum_{m=1}^{N_G}w_m=1,
\)
and $\mathcal N(\mu_m,\sigma_m^2)$ denotes a Gaussian distribution with mean $\mu_m$ and variance $\sigma_m^2$.
The model parameters
\begin{align}
   \left\{w_m,\mu_m,\sigma_m^2\right\}_{m=1}^{N_G}
\end{align}
are estimated from the training dataset $\mathcal D_k$ using the Expectation-Maximization (EM) algorithm, which iteratively maximizes the likelihood of the observed samples under the mixture model.
The number of mixture components governs the trade-off between model flexibility and complexity. To guide its selection, Gaussian mixtures with varying numbers of components were evaluated using the Bayesian Information Criterion (BIC). Across the considered fragmentation scenarios and propagation epochs, the optimal number of components was consistently found to lie in the range of four to six components. Consequently, a fixed value of
$N_G=5$
was adopted throughout this work \cite{Paul2025ICAWFragmentation,ParienteServitje2025CloudCollision}. Figure~\ref{fig:envisat_sma} illustrates the Gaussian Mixture Model approximation of the semi-major-axis distribution for the ENVISAT fragmentation cloud two years after breakup (see Fragmentation scenarios in \cref{sec:fragmentation_scenarios}). The empirical distribution exhibits a multimodal structure resulting from the combined effects of atmospheric drag and orbital perturbations on fragments with different area-to-mass ratios. The proposed GMM accurately captures both the dominant modes and the overall shape of the distribution using a limited number of mixture components, providing a compact parametric representation suitable for subsequent conditional modeling of the remaining orbital elements.

\subsection{Conditional Structure of the Cloud}

The principal limitation of the band-formation assumption discussed in Section~\ref{sec:introduction} is that it neglects potential dependencies between the angular variables and the remaining orbital elements. However, inspection of propagated fragmentation clouds reveals persistent structures linking the eccentricity, \gls{raan}, and \gls{aop} to the semi-major axis, even after long propagation times \cite{lorenzo2024space}.
To preserve the dominant dependencies induced by orbital dynamics, the conditional density
\begin{align}
p_k(e,\Omega,\omega~|~a)
\end{align}
is estimated from the propagated representative fragments rather than replaced by a uniform approximation. Direct estimation of this three-dimensional conditional density would however require a large number of parameters and a correspondingly large training dataset. To obtain a tractable surrogate model, we assume conditional independence given the semi-major axis and approximate 
\begin{align}
p_k(e,\Omega,\omega|a)
\approx
p_k(e|a)\,
p_k(\Omega|a)\,
p_k(\omega|a)
\label{eq:conditional_factorization}
\end{align}
This approximation should not be interpreted as a global independence assumption. On the contrary, dependencies between the orbital elements are preserved through their common dependence on the semi-major axis, which acts as the root variable of the hierarchical model. The resulting formulation significantly reduces model complexity while retaining the principal structures observed in the propagated debris clouds.
The following subsections describe the estimation procedures adopted for $p_k(e|a)$, $p_k(\Omega|a)$, and $p_k(\omega|a)$.

\subsection{Conditional Distribution of Eccentricity}
\label{sec:eccentricity}

\begin{figure}[H]
    \centering
    \begin{subfigure}[b]{0.48\textwidth}
        \centering
        \includegraphics[width=\linewidth]{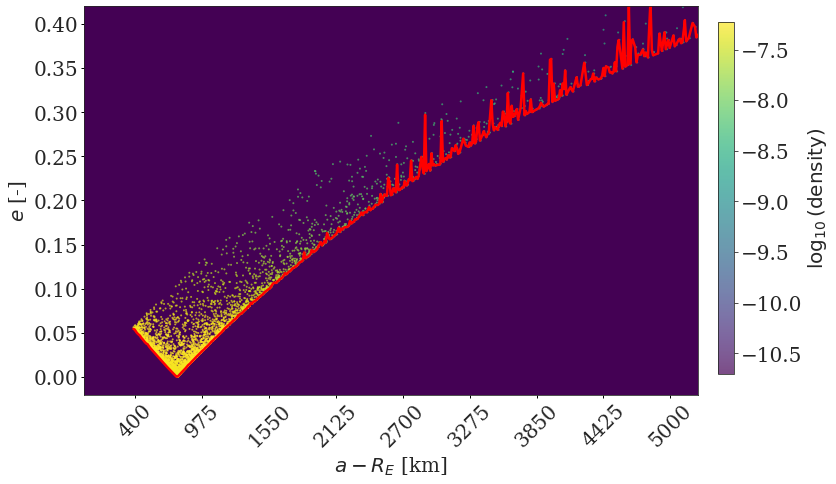}
        \caption{Fragmentation epoch.}
        \label{fig:envisat_e_wrt_sma_t0}
    \end{subfigure}
    \hfill
    \begin{subfigure}[b]{0.48\textwidth}
        \centering
        \includegraphics[width=\linewidth]{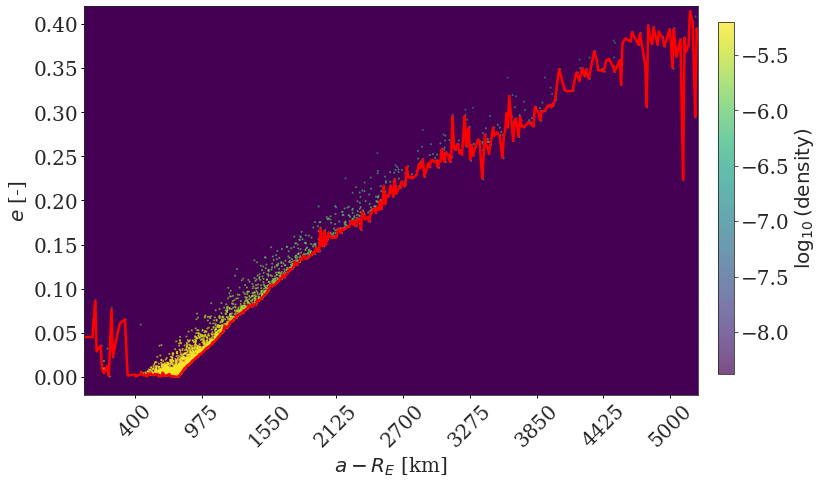}
        \caption{2 years after fragmentation.}
        \label{fig:envisat_e_wrt_sma_t2}
    \end{subfigure}
    \caption{\small $(a,e)$-distribution for the ENVISAT fragmentation scenario. The lower boundary is in red.}
    \label{fig:envisat_e_wrt_sma}
\end{figure}

The conditional distribution of eccentricity exhibits a strong dependence on the semi-major axis, as illustrated in Fig.~\ref{fig:envisat_e_wrt_sma}. Rather than being distributed over the entire $(a,e)$ plane, the propagated debris cloud is concentrated above a well-defined lower boundary whose shape vary slowly throughout the propagation. This observation suggests decomposing the conditional distribution into a deterministic component represented by the lower boundary and a residual stochastic component describing the dispersion around it.

Let
$e_{\min,k}(a)$
denote the estimated lower boundary at epoch $t_k$.
For the fragmentation epoch $t_0$, an analytical expression of
$e_{\min,0}(a)$
can be derived from the breakup dynamics \cite{giudici2023probabilistic}.
However, deriving a similar expression at later epochs would require accounting for the long-term evolution of the cloud under orbital perturbations. As illustrated in \cref{fig:envisat_e_wrt_sma_t2}, the resulting boundary exhibits a complex shape that is difficult to describe analytically.  Consequently, the lower boundary is approximated directly from the propagated samples using a piecewise-linear approximation. 
Since no analytical expression for this boundary 
is available, it is inferred directly from the 
propagated samples using a piecewise-linear 
approximation.
The eccentricity is then decomposed as
\begin{align}
e = e_{\min,k}(a) + r_e,
\label{eq:ecc_residual}
\end{align}
where $r_e \ge 0$ denotes a residual eccentricity. This transformation removes the dominant deterministic structure of the conditional distribution and yields a substantially simpler residual distribution \cite{Paul2025ICAWFragmentation}.

\begin{figure}[H]
    \centering
    \begin{subfigure}[b]{0.48\textwidth}
        \centering
        \includegraphics[width=\linewidth]{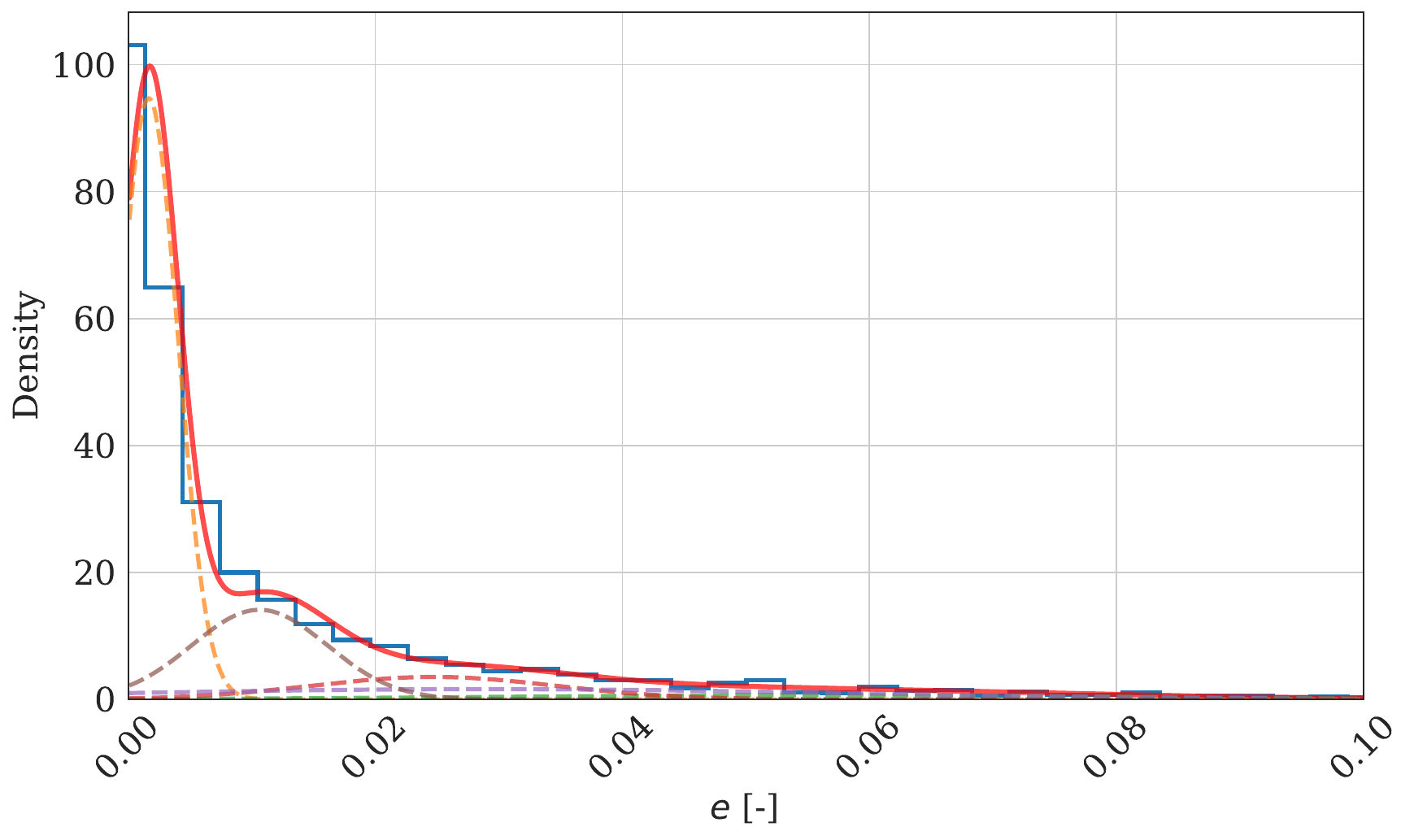}
        \caption{Fragmentation epoch.}
        \label{fig:envisat_shifted_e_wrt_sma_t0}
    \end{subfigure}
    \hfill
    \begin{subfigure}[b]{0.48\textwidth}
        \centering
        \includegraphics[width=\linewidth]{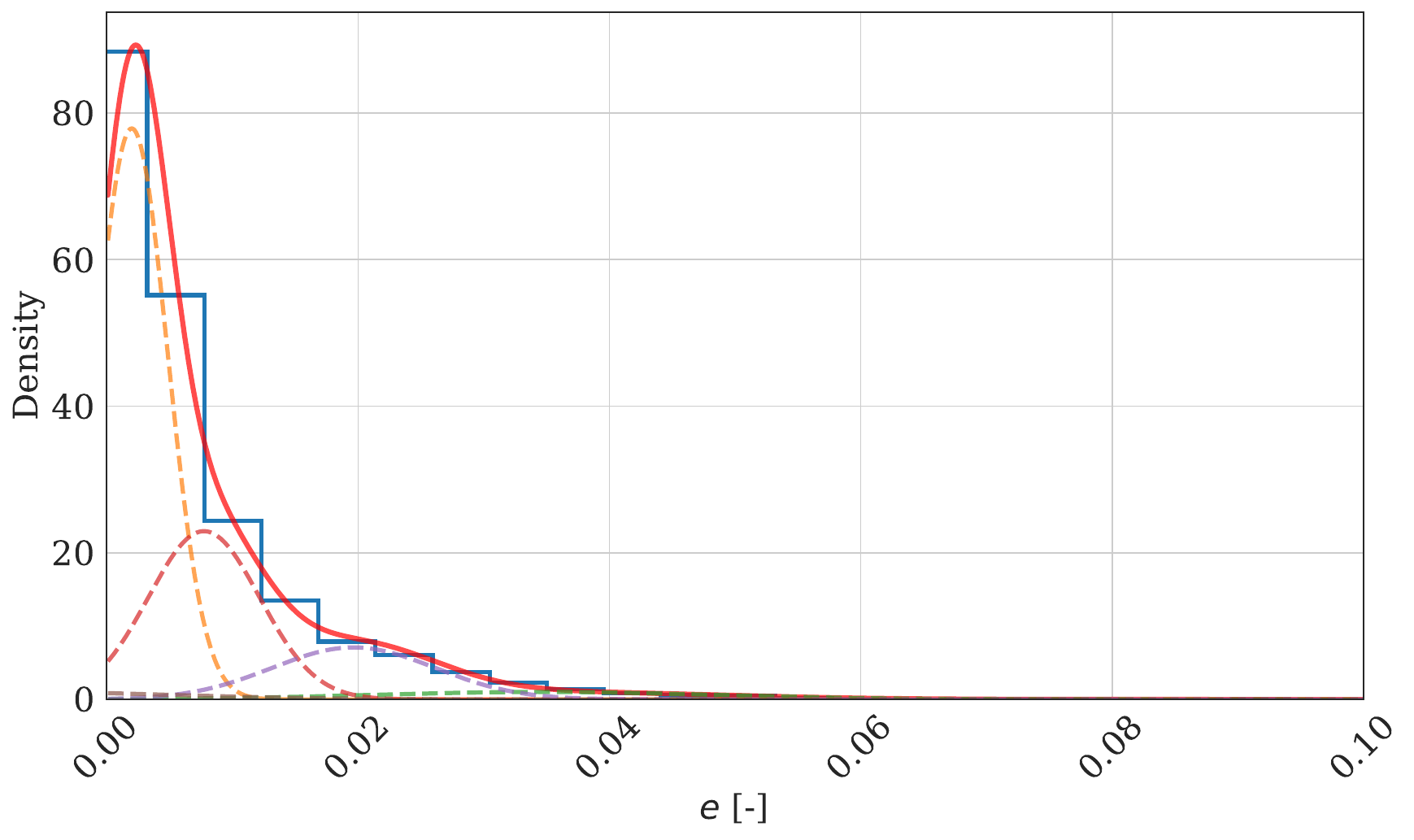}
        \caption{2 years after fragmentation.}
        \label{fig:envisat_shifted_e_wrt_sma_t2}
    \end{subfigure}

    \caption{\small Residual eccentricity distribution $r_e$ and corresponding Gaussian-mixture approximation.}
    \label{fig:envisat_shifted_e_wrt_sma}
\end{figure}

Once the lower boundary has been subtracted from the propagated samples, the residual variable $r_e$ is modelled using a Gaussian Mixture Model,
\begin{align}
p_k(r_e)=
\sum_{m=1}^{N_G}
w^{e}_{k,m}\,
\mathcal N\left(r_e;\mu^{e}_{k,m},(\sigma^{e}_{k,m})^2\right),
\label{eq:gmm_ecc_residual}
\end{align}
where the mixture weights, means, and variances are estimated by the EM algorithm. As for the semi-major-axis model, a fixed value $N_G=5$ is used throughout this work. Figure~\ref{fig:envisat_shifted_e_wrt_sma} illustrates the resulting fit. The empirical residual distribution is shown as a blue histogram, the overall Gaussian-mixture approximation as a red curve, and the individual mixture components as dashed lines. A small number of components is sufficient to capture the main features of the residual distribution across the considered propagation times.
The conditional density $p_k(e|a)$ is therefore obtained by shifting the estimated residual distribution according to the lower boundary,
\begin{align}
p_k(e|a)=
p_k\left(r_e=e-e_{\min,k}(a)\right).
\label{eq:ecc_conditional_shift}
\end{align}
Equivalently, sampling from $p_k(e|a)$ is performed by first evaluating $e_{\min,k}(a)$ and then drawing a residual eccentricity $r_e$ from the fitted mixture model.

\subsection{Conditional Distribution of RAAN}
\label{subsec:raan}

\begin{figure}[H]
    \centering
    \begin{subfigure}[b]{0.48\textwidth}
        \centering
        \includegraphics[width=\linewidth]{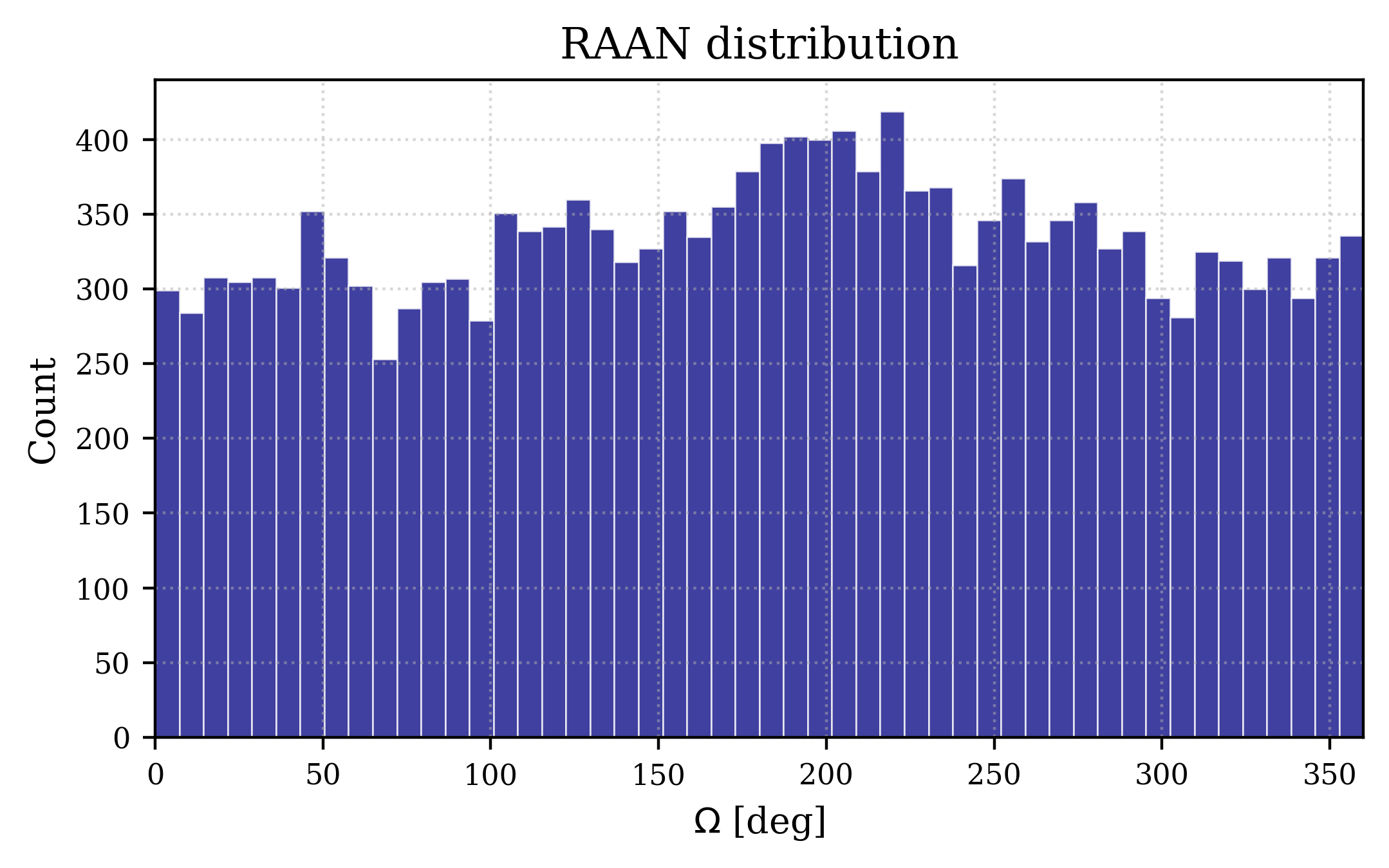}
        \caption{Marginal $\Om$-distribution.}
        \label{fig:raan_hist_2y}
    \end{subfigure}
    \hfill
    \begin{subfigure}[b]{0.48\textwidth}
        \centering
        \includegraphics[width=\linewidth]{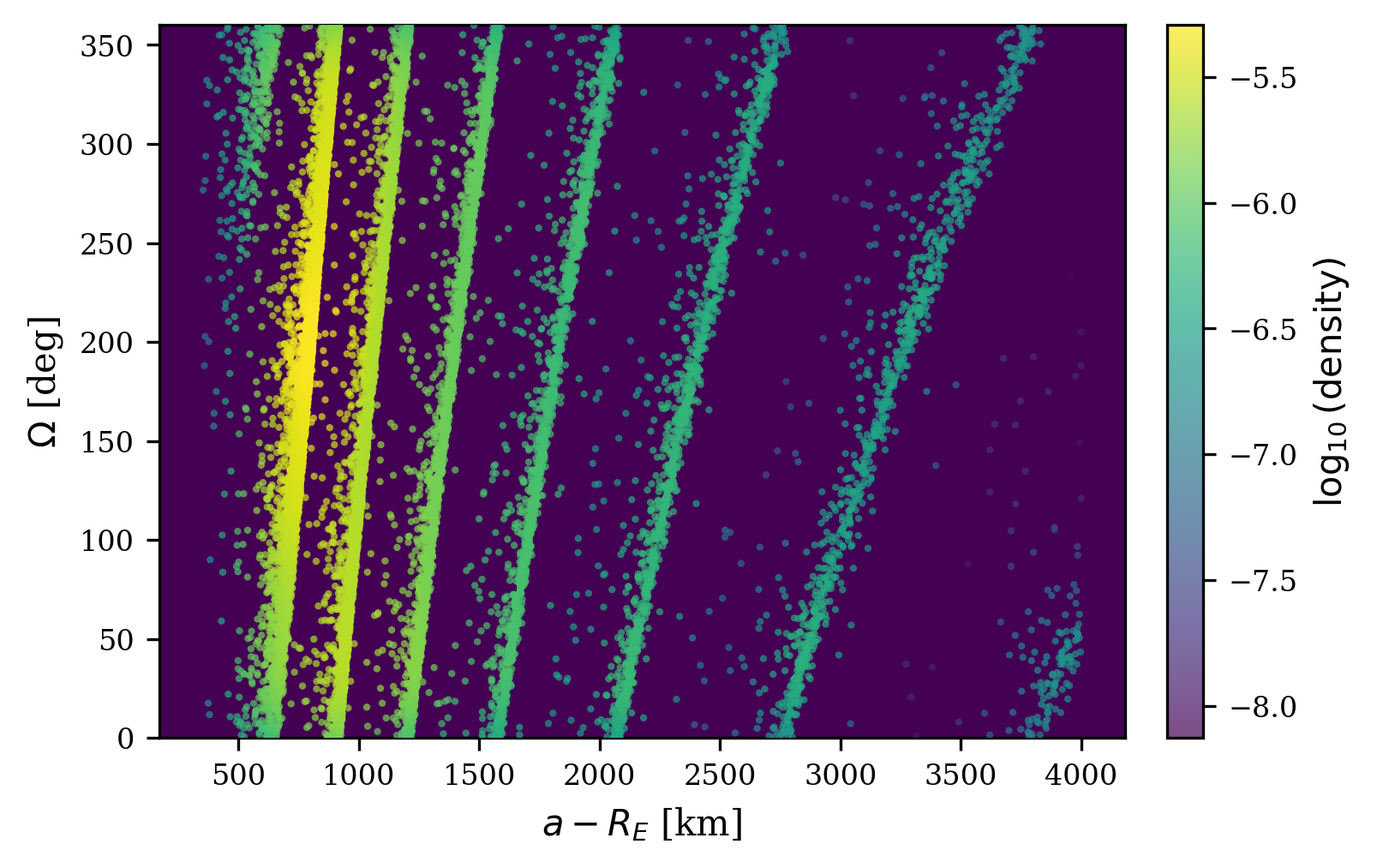}
        \caption{$(a,\Om)$-distribution.}
        \label{fig:raan_a_2y}
    \end{subfigure}

    \caption{Distribution of the \gls{raan}~$\Omega$ 
    after 2~years for the \textbf{LeoFrag} fragmentation scenario. }
    \label{fig:raan_2y}
\end{figure}

\begin{figure}[H]
    \centering
    \begin{subfigure}[b]{0.48\textwidth}
        \centering
        \includegraphics[width=\linewidth]{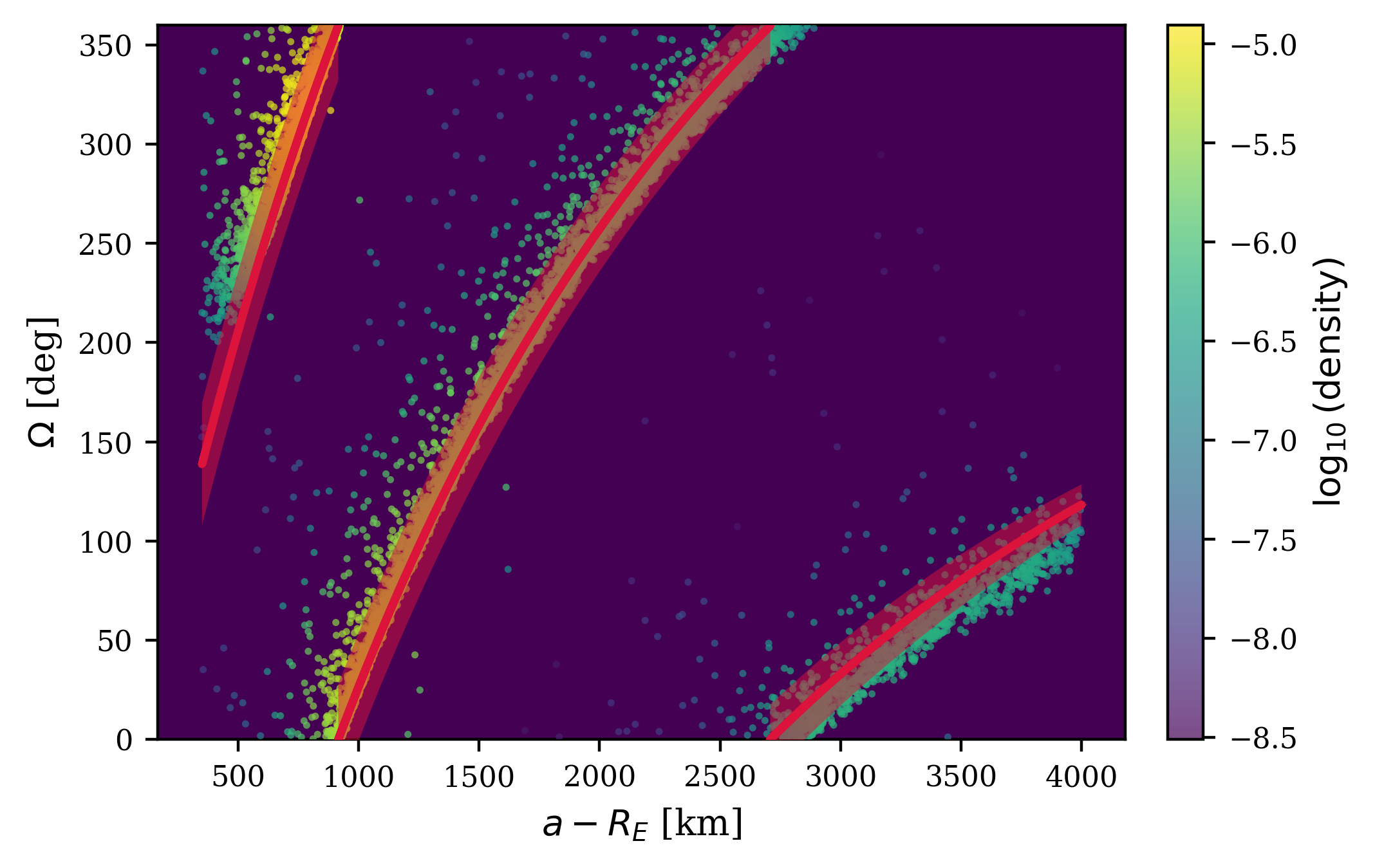}
        \caption{6 months after fragmentation.}
        \label{fig:leo_raan_wrt_sma_t0}
    \end{subfigure}
    \hfill
    \begin{subfigure}[b]{0.48\textwidth}
        \centering
        \includegraphics[width=\linewidth]{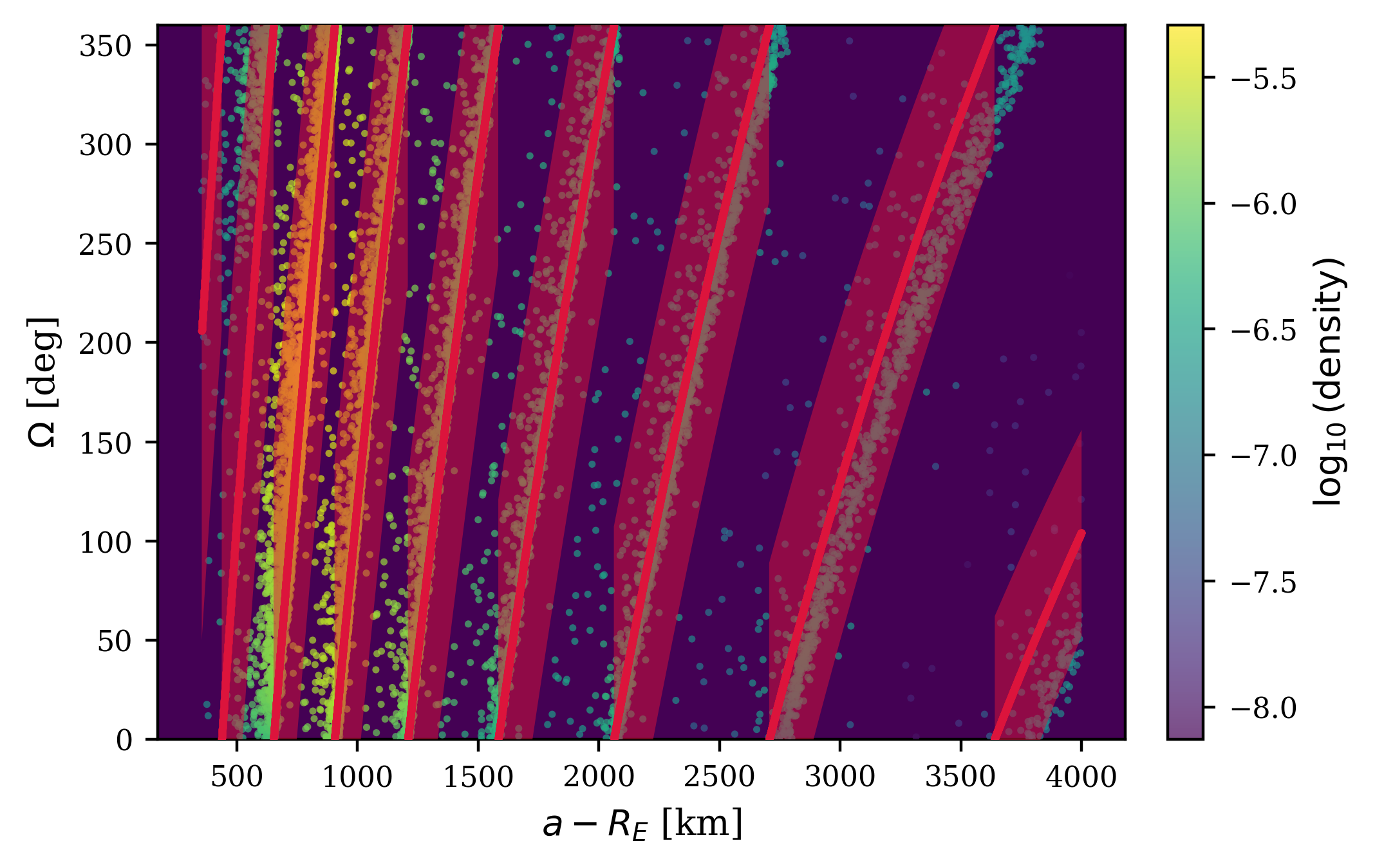}
        \caption{2 years after fragmentation.}
        \label{fig:leo_raan_wrt_sma_t2}
    \end{subfigure}

    \caption{Fitting of $(a,\Om)$-distribution for the \textbf{LeoFrag} fragmentation scenario with~2$\sigma$ confidence envelope.}
    \label{fig:leo_raan_wrt_sma}
\end{figure}

The conditional distribution of the \gls{raan} exhibits a markedly different structure from that of eccentricity. 
Figure~\ref{fig:raan_2y} illustrates the motivation for
introducing a conditional model for the RAAN.
Although the marginal distribution of $\Omega$
appears nearly uniform after two years of propagation
(Fig.~\ref{fig:raan_hist_2y}),
strong structures remain visible when conditioning
on the semi-major axis (Fig.~\ref{fig:raan_a_2y}).
These structures arise from differential nodal precession
and cannot be captured by models assuming independence
between $a$ and $\Omega$.
As shown in Fig.~\ref{fig:leo_raan_wrt_sma}, the debris population is concentrated around narrow curved bands in the $(a,\Omega)$ plane. These structures remain clearly visible after several months and even years of propagation, demonstrating that the semi-major axis retains a strong influence on the RAAN distribution.

The observed bands correspond to continuous curves once the angular variable is unwrapped. For a generic angular variable $\theta$, we define
\begin{equation}
\widetilde{\theta}_{t}
=
\theta_{t}
+
2\pi n_t,
\qquad
n_t\in\mathbb Z,
\label{eq_angle_unwrap}
\end{equation}
where $n_t$ is chosen such that the mapping
$t\mapsto \widetilde{\theta}_t$
remains continuous.
After unwrapping, the conditional distribution exhibits a clear functional relationship between the semi-major axis and the RAAN. The conditional mean is approximated by
\begin{equation}
\mu^{\Omega}_{k}(a)
=
\frac{c^{\Omega}_{k}}
{\left(a-a^{\Omega}_{k}\right)^{d^{\Omega}_{k}}},
\label{eq_mean_Omega}
\end{equation}
where the parameters
$c^{\Omega}_{k}$,
$a^{\Omega}_{k}$,
and
$d^{\Omega}_{k}$
are estimated from the propagated samples using nonlinear least-squares regression \cite{Paul2025FragmentationReport,Bazenga2024BandFormation}.
The remaining dispersion around the fitted curve is modeled as a Gaussian residual,
\begin{align}
p_k(\Omega|a)
=
\mathcal N
\!\left(
\Omega;
\mu_k^\Omega(a),
\sigma_{\Omega,k}^2
\right),
\label{eq_density_raan_wrt_sma}
\end{align}
where $\sigma_{\Omega,k}^2$ is estimated from the regression residuals after Z-score-based outlier rejection \cite{ParienteServitje2025CloudCollision}.
Sampling from $p_k(\Omega|a)$ is then performed by evaluating $\mu_k^\Omega(a)$ and adding a Gaussian residual drawn according to Eq.~\eqref{eq_density_raan_wrt_sma}.

\subsection{Conditional Distribution of \gls{aop}}
\label{subsec:aop}

\begin{figure}[H]
    \centering
    \begin{subfigure}[b]{0.48\textwidth}
        \centering
        \includegraphics[width=\linewidth]{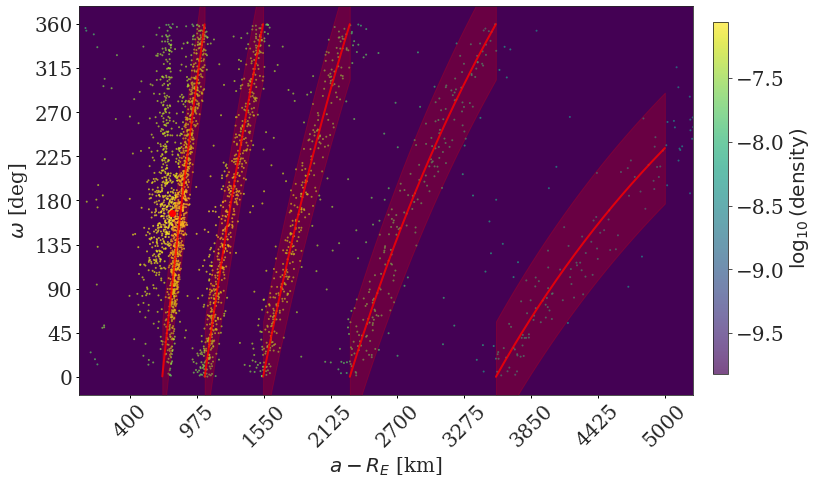}
        \caption{\small 2 years after fragmentation.}
        \label{fig:envisat_aop_wrt_sma_t2}
    \end{subfigure}
    \hfill
    \begin{subfigure}[b]{0.48\textwidth}
        \centering
        \includegraphics[width=\linewidth]{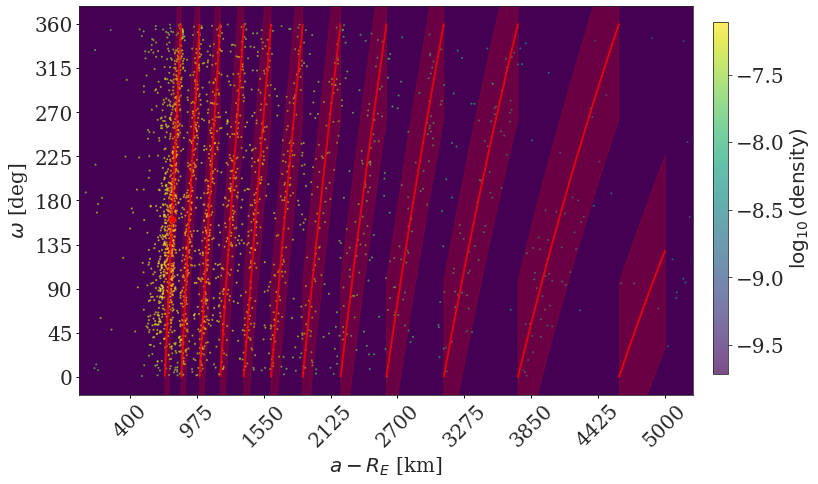}
        \caption{\small 5 years after fragmentation.}
        \label{fig:envisat_aop_wrt_sma_t5}
    \end{subfigure}
    \caption{\small $(a,\om)$-distribution for the ENVISAT fragmentation scenario. The fitted mean curve and $2\sigma$ envelope are shown in red.}
    \label{fig:envisat_aop_wrt_sma}
\end{figure}

The conditional distribution of the \gls{aop} exhibits a structure similar to that of the \gls{raan}, with debris concentrated around narrow curved bands in the $(a,\omega)$ plane. Examples are shown in Fig.~\ref{fig:envisat_aop_wrt_sma}. 
As for $\Om$, these patterns reveal a strong dependence between $\omega$ and $a$ that is not captured by uniform-band approximations.
After angle unwrapping using Eq.~(\ref{eq_angle_unwrap}), the upper branch of the conditional distribution is represented by the parametric model
\begin{equation}
\mu_{k}^{\omega}(a)
=
\frac{c_{k}^{\omega}}
{\left(a-a_{k}^{\omega}\right)^{d_{k}^{\omega}}}
+s_k^\omega,
\label{eq_mean_omega}
\end{equation}
where the parameters
$c_k^\omega$,
$a_k^\omega$,
$d_k^\omega$,
and
$s_k^\omega$
are estimated from the propagated samples using nonlinear least-squares regression.
The remaining dispersion around the fitted curve is modeled as a Gaussian residual,
\begin{align}
p_k^{>}(\omega|a)
=
\mathcal N
\!\left(
\omega;
\mu_k^\omega(a),
\sigma_{\omega,k}^2
\right),
\label{eq_density_aop_wrt_sma_upper}
\end{align}
where $\sigma_{\omega,k}^2$ is estimated from the regression residuals after Z-score-based outlier rejection.

\begin{figure}[H]
    \centering
    \includegraphics[width=0.6\linewidth]{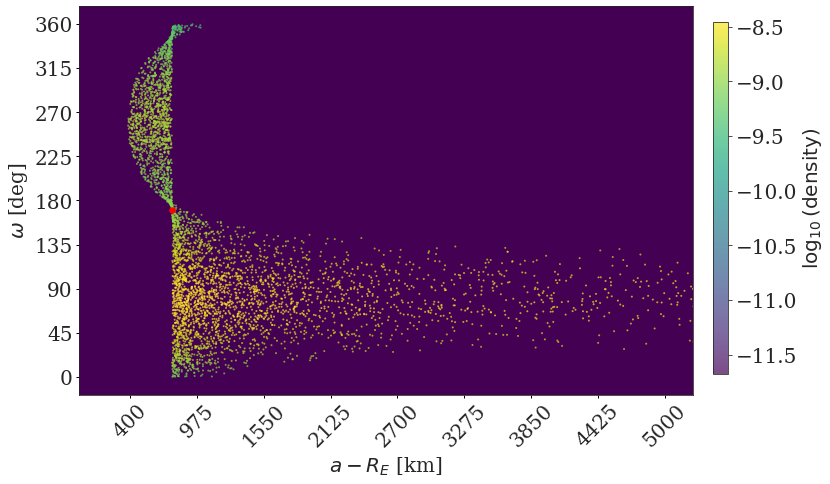}
\caption{\small $(a,\om)$-distribution at fragmentation epoch. The symmetry center separating the upper and lower branches is shown in red.}
\label{fig:envisat_aop_split}
\end{figure}

A closer inspection of the cloud reveals that a single branch is insufficient to represent the full conditional distribution. As illustrated in Fig.~\ref{fig:envisat_aop_split}, the debris population naturally separates into two branches around a symmetry center associated with the fragmenting object. The upper branch corresponds to debris with semi-major axis greater than that of the fragmenting object, whereas the lower branch corresponds to debris with smaller semi-major axis. This separation originates from fragments receiving positive or negative orbital-energy variations during the breakup, leading respectively to semi-major axes larger or smaller than that of the fragmenting object.

The upper branch corresponds to the structure modeled by Eq.~(\ref{eq_density_aop_wrt_sma_upper}). To obtain a simple representation of the lower branch, the model exploits the approximate symmetry of the cloud around the symmetry center, depicted in red in Fig.~\ref{fig:envisat_aop_split}. Empirically, this point was found to closely follow the nominal trajectory $(a_t^{\mathrm p},\omega_t^{\mathrm p})$ of the parent object, i.e., the trajectory that would have been obtained had the fragmentation not occurred at epoch $t_0$. Let $(a_k^{\mathrm p},\omega_k^{\mathrm p})$ denote the location of this point at epoch $t_k$.

The lower branch is then obtained through the symmetry transformation
\begin{equation}
p_k^{<}(\omega|a)
=
p_k^{>}
\!\left(
2\omega_k^{\mathrm p}-\omega
\mid
2a_k^{\mathrm p}-a
\right).
\end{equation}
The resulting model provides a compact representation of both branches 
while preserving the principal geometric features observed in the 
propagated cloud \cite{Paul2025ICAWFragmentation}. Although the symmetry assumption is only approximate, it was found sufficient to reproduce the principal geometric features of the fragmentation clouds considered in this work.

\subsection{Inclination and Mean Anomaly}

Consistent with most statistical fragmentation-cloud models in the literature, the inclination is assumed constant over the propagation horizons considered in this work \cite{letizia_space_2016,Frey_S_2021_1,paul24fragmentation}. Consequently, it is treated as independent of the remaining orbital elements and modeled through the Dirac distribution
\begin{equation}
p_k(i)
\coloneq
\delta_{i_0^{\mathrm p}}(i),
\end{equation}
where $i_0^{\mathrm p}$ denotes the inclination of the parent object at fragmentation epoch $t_0$.

Similarly, the mean anomaly is assumed independent of the remaining orbital elements and uniformly distributed over $[0,2\pi)$, following the standard phased-cloud description adopted in the literature \cite{letizia_space_2016,Frey_S_2021_1,paul24fragmentation}. This assumption is motivated by the rapid along-track dispersion of the debris cloud after fragmentation. Differential orbital periods quickly spread the fragments around the orbit, causing the mean anomaly distribution to approach uniformity within a few days, a timescale negligible compared with the multi-year horizons considered in this work. Accordingly,
\begin{equation}
p_k(M)
=
\frac{1}{2\pi},
\qquad
M\in[0,2\pi).
\end{equation}
Consequently, the proposed model focuses its complexity on the orbital elements that exhibit strong conditional structures, namely the eccentricity, RAAN, and argument of perigee.

\subsection{Generation of Synthetic Debris Clouds}

Once the distributions introduced in the previous subsections have been fitted at epoch $t_k$, synthetic debris states can be generated without additional orbital propagation.
The generation procedure follows the hierarchical factorization
\begin{equation}
p_k(a,e,i,\Omega,\omega,M)
=
p_k(a)\,
p_k(e|a)\,
p_k(\Omega|a)\,
p_k(\omega|a)\,
p_k(i)\,
p_k(M),
\label{eq:hgdm_factorization}
\end{equation}
and a synthetic debris state $\tilde x=(\tilde a,\tilde e,\tilde i,\tilde \Om,\tilde \om,\tilde M)$ is generated according to the following procedure:
\begin{enumerate}
\item Sample the semi-major axis
$\tilde a \sim p_k(a).$
\item Sample, conditionally on $\{a=\tilde a\}$, the eccentricity; RAAN and AOP such that 
\item Sample, conditionally on $a=\tilde a$,
\begin{align}
\tilde e &\sim p_k(e\,|\,a=\tilde a),&
\tilde\Omega &\sim p_k(\Omega\,|\,a=\tilde a),&\text{and   }
\tilde\omega \sim p_k(\omega\,|\,a=\tilde a).
\end{align}
\item Assign the inclination according to and sample the mean anomaly according to
$\tilde M \sim \mathcal U([0,2\pi))$, the uniform distribution on $[0,2\pi).$
\end{enumerate}
Repeating this procedure yields an arbitrary number of synthetic debris states distributed according to the surrogate density $p_k$. The resulting cloud can subsequently be used for density estimation, collision-risk assessment, and long-term population studies without requiring additional high-fidelity propagations.

The present work focuses on the spatial distribution of the debris cloud. Population-size estimation can be incorporated separately using the point-process framework developed in \cite{paul24fragmentation}. Given an initial population size $N_0^{\mathrm{deb}}$ provided by the NASA Breakup Model, a first-order estimate of the surviving population at epoch $t_k$ is
\begin{equation}
N_k^{\mathrm{deb}}
=
N_0^{\mathrm{deb}}
\int_{\{h_{\mathrm p}>h_{\mathrm{decay}}\}}
p_k(x)\,dx.
\end{equation}
where $R_{\mathrm{E}}$ represents Earth radius, 
$h_{\p}\coloneq a(1-e)-R_{E}$ the perigee altitude of a debris and 
$h_{\mathrm{decay}}$ an altitude where we consider the debis as decayed, ( e.g for $h_{\mathrm{decay}}=150$ km).

\section{Evaluation Setup}
\label{sec:evaluation_setup}

This section describes the fragmentation scenarios, reference distributions, and evaluation protocol used to assess the proposed Hierarchical Generative Density Model (HGDM). The objective is to compare the estimated distributions $p_k$ with high-fidelity reference distributions $p_k^{\mathrm{ref}}$ obtained from large-scale Monte Carlo simulations.

\subsection{Fragmentation Scenarios}
\label{sec:fragmentation_scenarios}

Three fragmentation scenarios are considered: representative breakup events based on ENVISAT and Ariane~2~R/B, and a synthetic low-Earth-orbit case, denoted \emph{LeoFrag}, introduced to isolate the secular mechanisms driving the dependence between semi-major axis and RAAN. LeoFrag corresponds to a circular orbit at 800 km altitude and $80^\circ$ inclination. The initial conditions of the three parent objects are summarized in Table~\ref{tab_simu}.
\begin{table}[h]
\centering
\caption{\small Initial properties of the fragmentation scenarios.}
\label{tab_simu}
\begin{tabular}{lcccc}
\toprule
Scenario & Mass [kg] & $a$ [km] & $e$ & $i$ [$^\circ$] \\
\midrule
ENVISAT      & 8110  & $R_e+765$ & 0.0001 & 98.3 \\
Ariane 2 R/B & 10000 & 24486     & 0.7188 & 8.4 \\
LeoFrag      & 2500  & $R_e+800$ & 0.0010 & 80.0 \\
\bottomrule
\end{tabular}
\end{table}
Only fragments with characteristic lengths between $\SI{1}{\centi\meter}$ and $\SI{1}{\meter}$ are considered, as this size range contains most debris relevant to long-term collision risk assessment while remaining largely untrackable \cite{giudici2024density,paul24fragmentation}. According to the NASA Breakup Model, this yields approximately 240,000, 180,000, and 100,000 fragments for ENVISAT, Ariane~2~R/B, and LeoFrag, respectively.
Each cloud is propagated for 5 years and sampled every 30 days. Reference trajectories are generated with the \gls{cnes} semi-analytical propagator STELA \cite{LeFevre_C_2012_1}, including atmospheric drag, solar radiation pressure, luni-solar perturbations, geopotential zonal harmonics ($J_2$--$J_4$), and tesseral resonances up to degree and order four.

\subsection{Reference Distributions}

The reference distributions $p_k^{\mathrm{ref}}$ are obtained through large-scale Monte Carlo simulations. For each fragmentation scenario, 200 independent Monte Carlo realizations are generated from the NASA Breakup Model and propagated using STELA.
These simulations provide high-fidelity estimates of the debris-cloud distributions at the successive evaluation epochs. The resulting Monte Carlo populations are regarded throughout this work as the ground-truth distributions against which the proposed surrogate model is evaluated.
Unlike the proposed HGDM model, which aims at estimating a compact probabilistic representation of the cloud, the reference approach explicitly propagates every simulated fragment. Consequently, it provides a highly accurate description of the cloud evolution at the expense of substantial computational and storage requirements.

\subsection{Training Procedure}

The proposed HGDM is trained using a limited subset of representative fragments propagated with the high-fidelity reference model. For each fragmentation scenario, an initial set of $N_{\mathrm{train}}=2000$ representative fragments is generated from the NASA Breakup Model and propagated to the successive evaluation epochs using STELA.

At each epoch $t_k$, the propagated fragments constitute the training dataset
$
\mathcal D_k
=
\left\{
{x}_k^{(j)}
\right\}_{j=1}^{N_{\mathrm{train}}}$,
where ${x}_k^{(j)}$ denotes the Keplerian state of fragment $j$ at epoch $t_k$.
The HGDM parameters are estimated independently at each evaluation epoch. The semi-major-axis and residual-eccentricity distributions are fitted using Gaussian Mixture Models with five components, while the conditional distributions of the angular variables are estimated through nonlinear regression followed by Gaussian residual modeling.
Once fitted, the resulting model defines a compact generative representation of the fragmentation cloud. Synthetic debris populations of arbitrary size can subsequently be generated without additional orbital propagation by sampling from the estimated probability distributions.
The reference Monte Carlo populations used for evaluation are not used during model fitting. Consequently, the assessment measures the ability of the surrogate model to reconstruct the statistical structure of large propagated fragmentation clouds from a limited number of representative samples.

\subsection{Evaluation Metrics}

The quality of the proposed surrogate model is assessed using both qualitative and quantitative criteria.
Qualitative evaluation is performed through pairwise projections of the orbital-element distributions. These visualizations provide insight into the ability of the estimated model to reproduce the multidimensional structure of the propagated fragmentation cloud.
Quantitative evaluation is performed using distributional distances computed between the reference density $p_k^{\mathrm{ref}}$ and the surrogate density $p_k$. Following previous work, the first-order error ($E_1$), second-order error ($E_2$), and symmetric Kullback--Leibler divergence ($D_{\mathrm{KL}}^{\mathrm{sym}}$) are computed on discretized density estimates.

To assess the benefit of modeling conditional dependencies, the proposed HGDM is also compared against a simplified baseline model in which the orbital elements are assumed independent and the angular variables are uniformly distributed, consistent with the classical band-formation approximation frequently adopted in long-term debris-cloud models.

\section{Results}
\label{sec:results}

This section evaluates the ability of the proposed Hierarchical Generative Density Model (HGDM) to reconstruct the long-term distribution of orbital fragmentation clouds from a limited number of propagated representative fragments. The assessment is performed through qualitative comparisons of the reconstructed distributions, quantitative distributional metrics, and measurements of computational and storage efficiency.

\subsection{Qualitative Reconstruction of Fragmentation Clouds}

     \begin{figure}[H]
    \centering
    \begin{subfigure}[b]{0.48\textwidth}
        \centering
        \includegraphics[width=\textwidth]{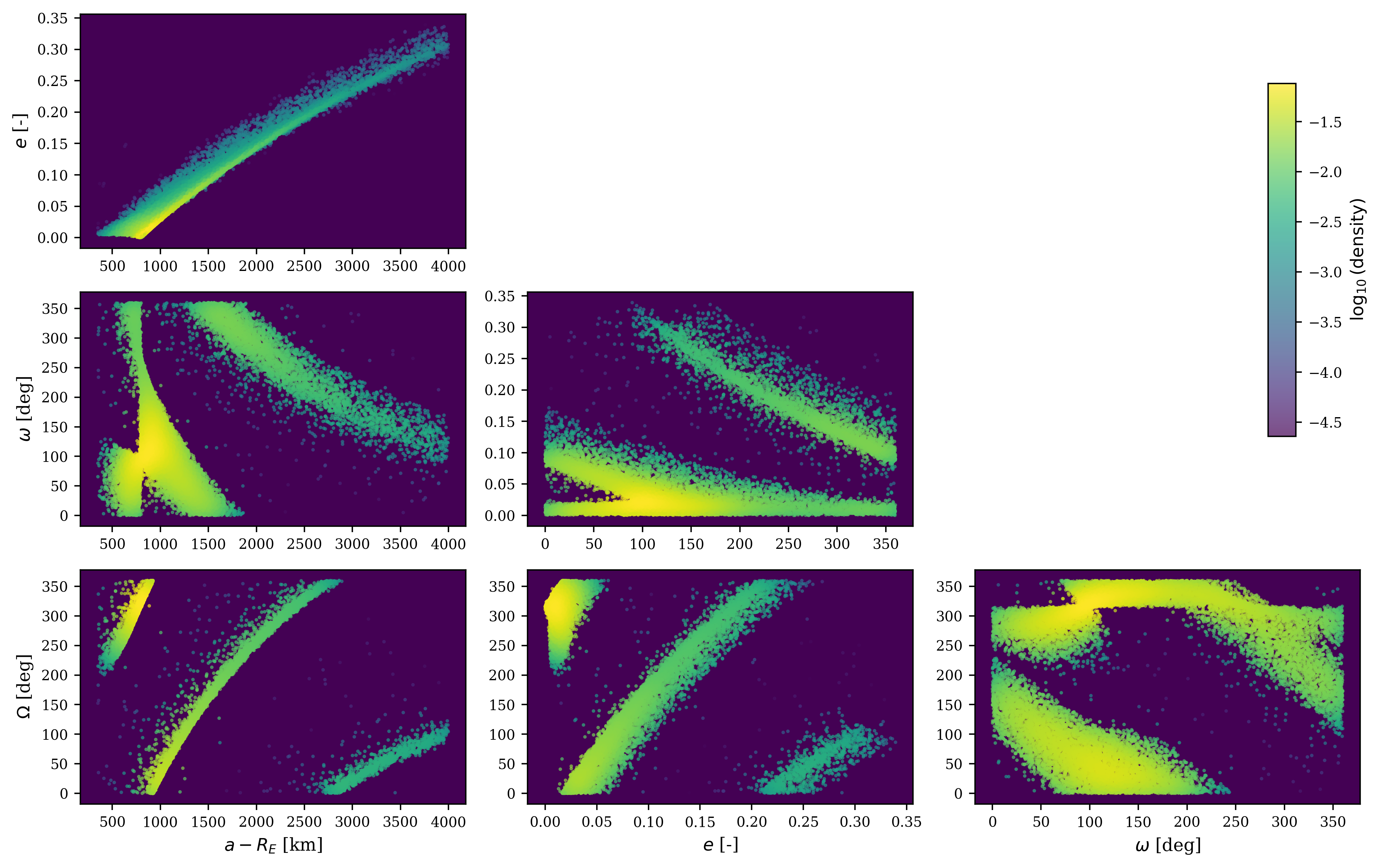}
        \caption{Reference -- 6 months.}
        \label{fig:leo_6m_ref}
    \end{subfigure}
    \hfill
    \begin{subfigure}[b]{0.48\textwidth}
        \centering
        \includegraphics[width=\textwidth]{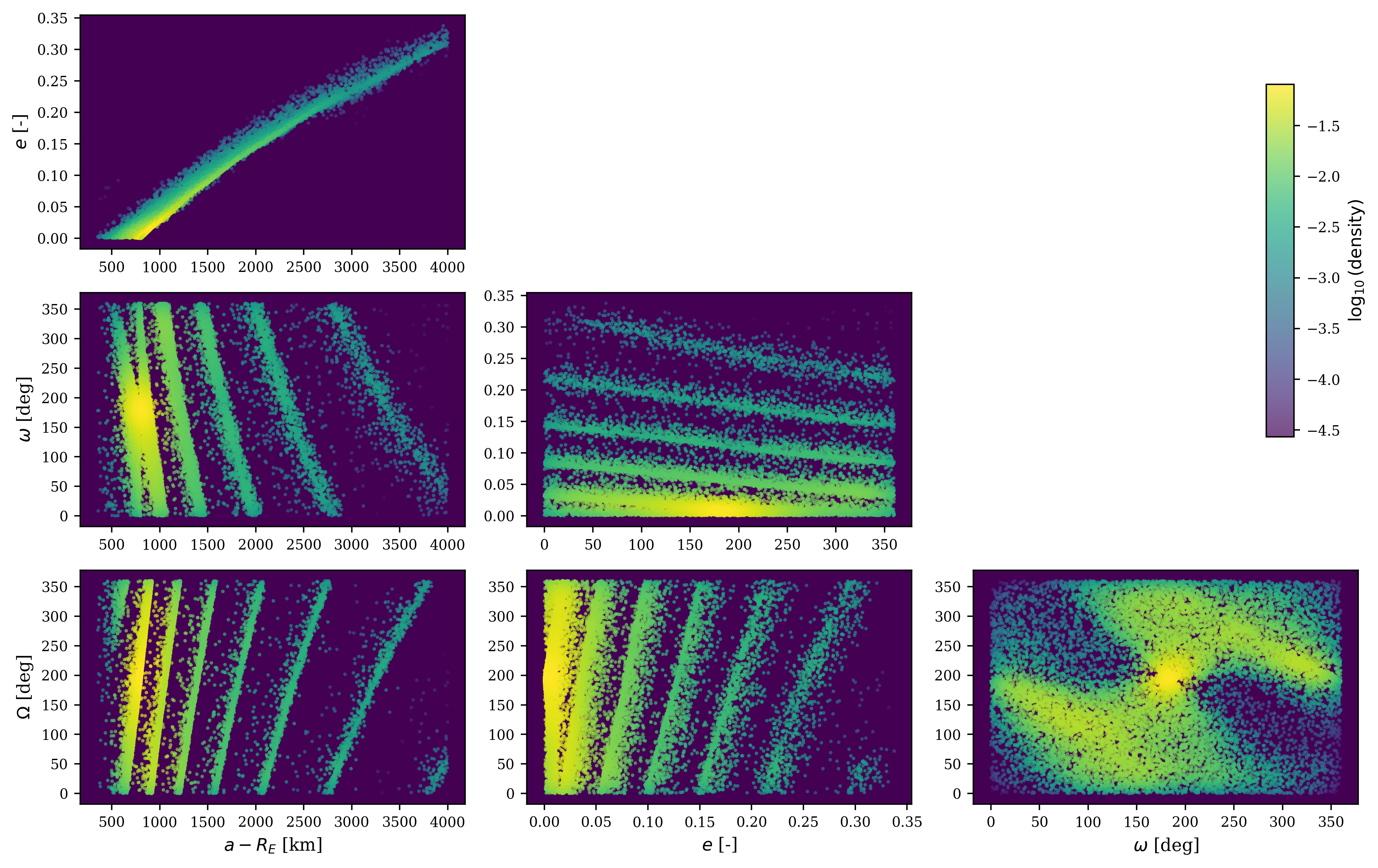}
        \caption{Reference -- 2 years.}
        \label{fig:leo_2y_ref}
    \end{subfigure}
    \\[0.5em] 
    \begin{subfigure}[b]{0.48\textwidth}
        \centering
        \includegraphics[width=\textwidth]{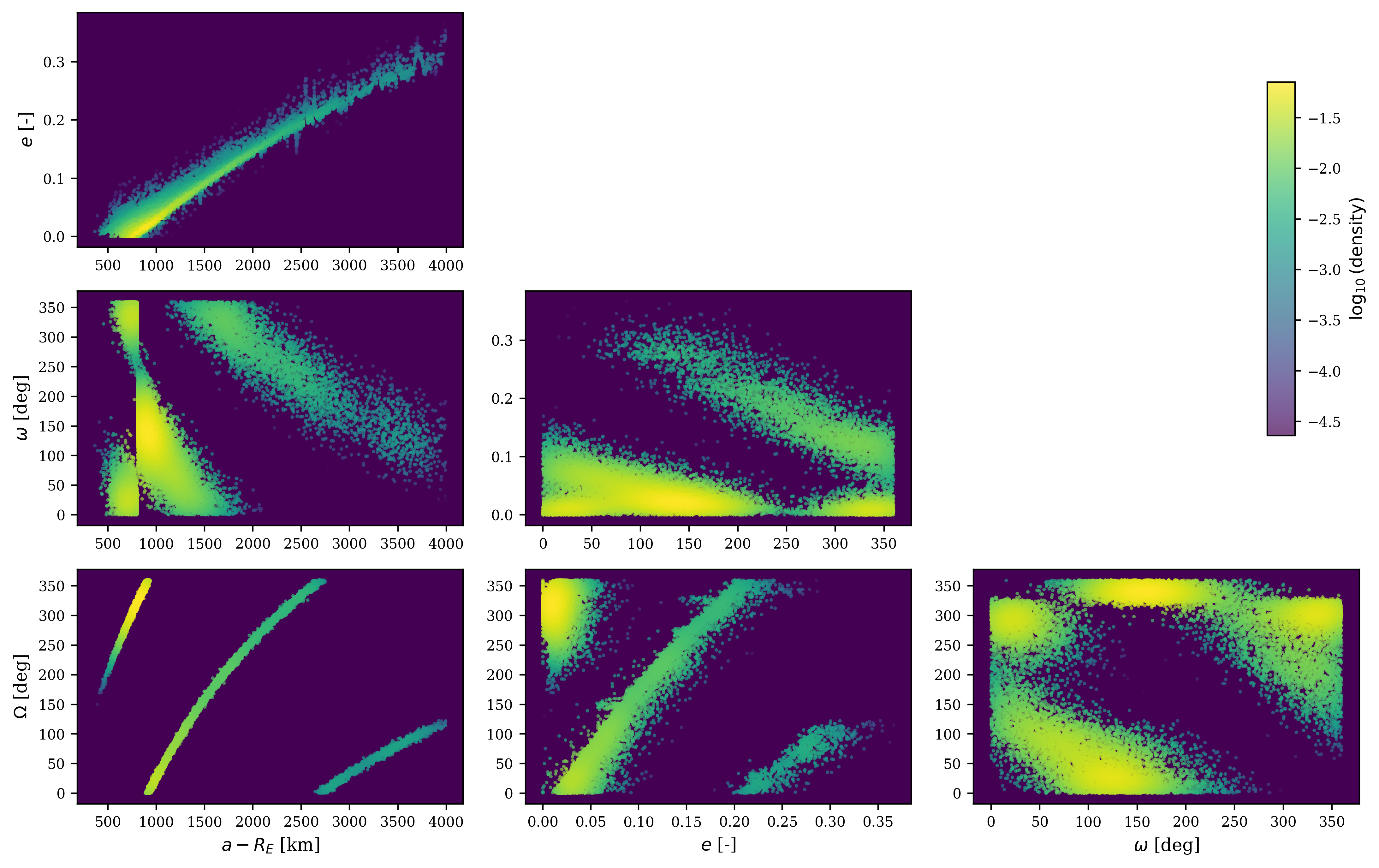}
        \caption{Surrogate -- 6 months.}
        \label{fig:leo_6m_pp}
    \end{subfigure}
    \hfill
    \begin{subfigure}[b]{0.48\textwidth}
        \centering
        \includegraphics[width=\textwidth]{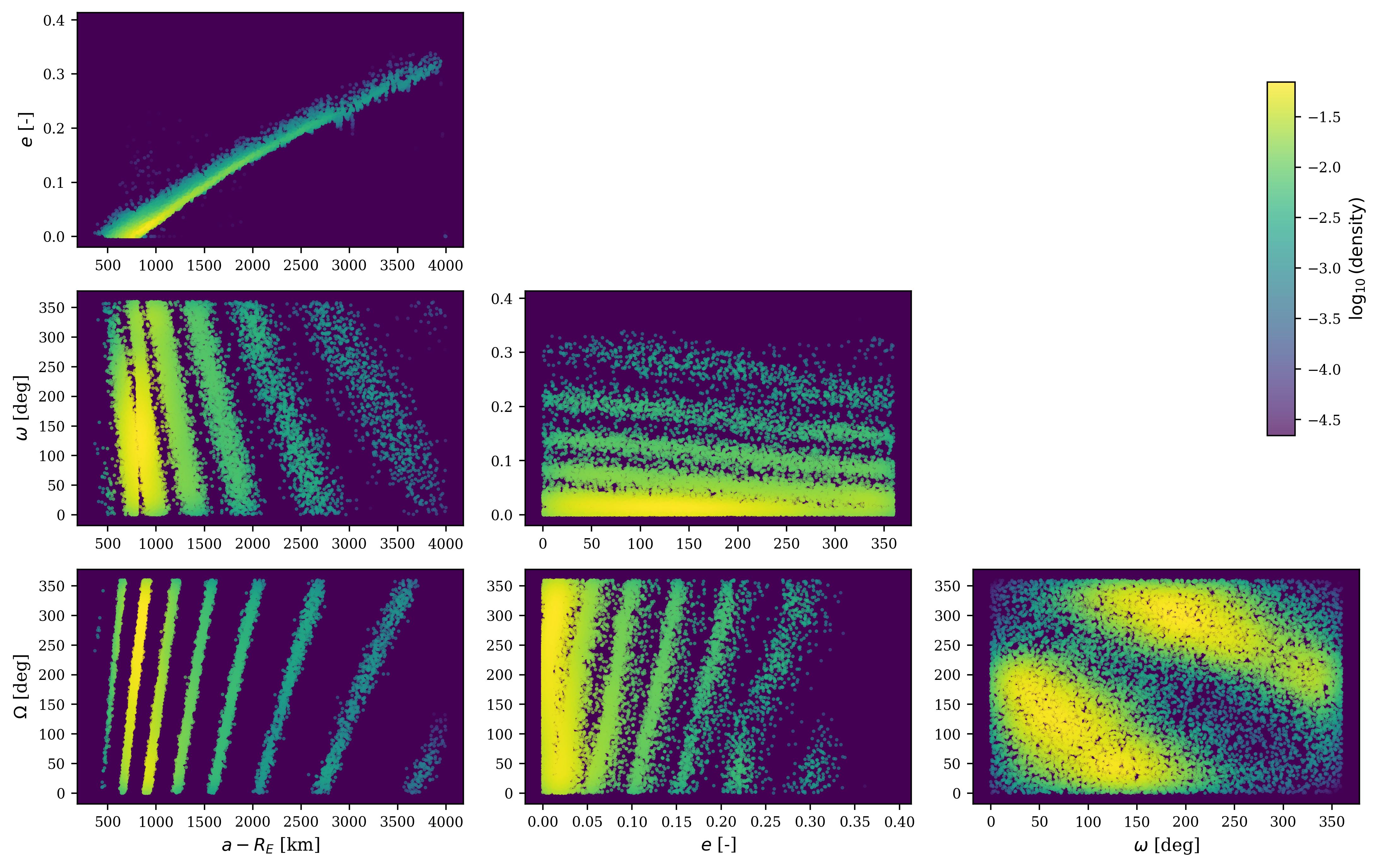}
        \caption{Surrogate -- 2 years.}
        \label{fig:leo_2y_pp}
    \end{subfigure}

    \caption{\small LeoFrag scenario. Pairwise orbital-element projections of the Monte Carlo reference (top) and the proposed surrogate (bottom).
.}
    \label{fig:leo_compact}
\end{figure}   

Figure~\ref{fig:leo_compact} compares the distributions generated by the proposed HGDM with high-fidelity Monte Carlo reference distributions for the LeoFrag scenario. The surrogate accurately reproduces the principal structures of the propagated cloud, including the lower-bound geometry of the $(a,e)$ distribution and the characteristic bands induced by differential precession in the $(a,\Omega)$ and $(a,\omega)$ planes \cite{Paul2025ICAWFragmentation}.
Interestingly, several correlations emerge naturally despite not being explicitly modeled. This suggests that the semi-major axis acts as an effective latent variable governing much of the long-term morphology of the fragmentation cloud.
Additional results for the Ariane~2~R/B scenario are provided in Appendix~\ref{app:ariane}.

\subsection{Distributional Accuracy}

\begin{figure*}[t]
\centering

\begin{subfigure}[t]{0.48\textwidth}
\centering
\includegraphics[width=\linewidth]
{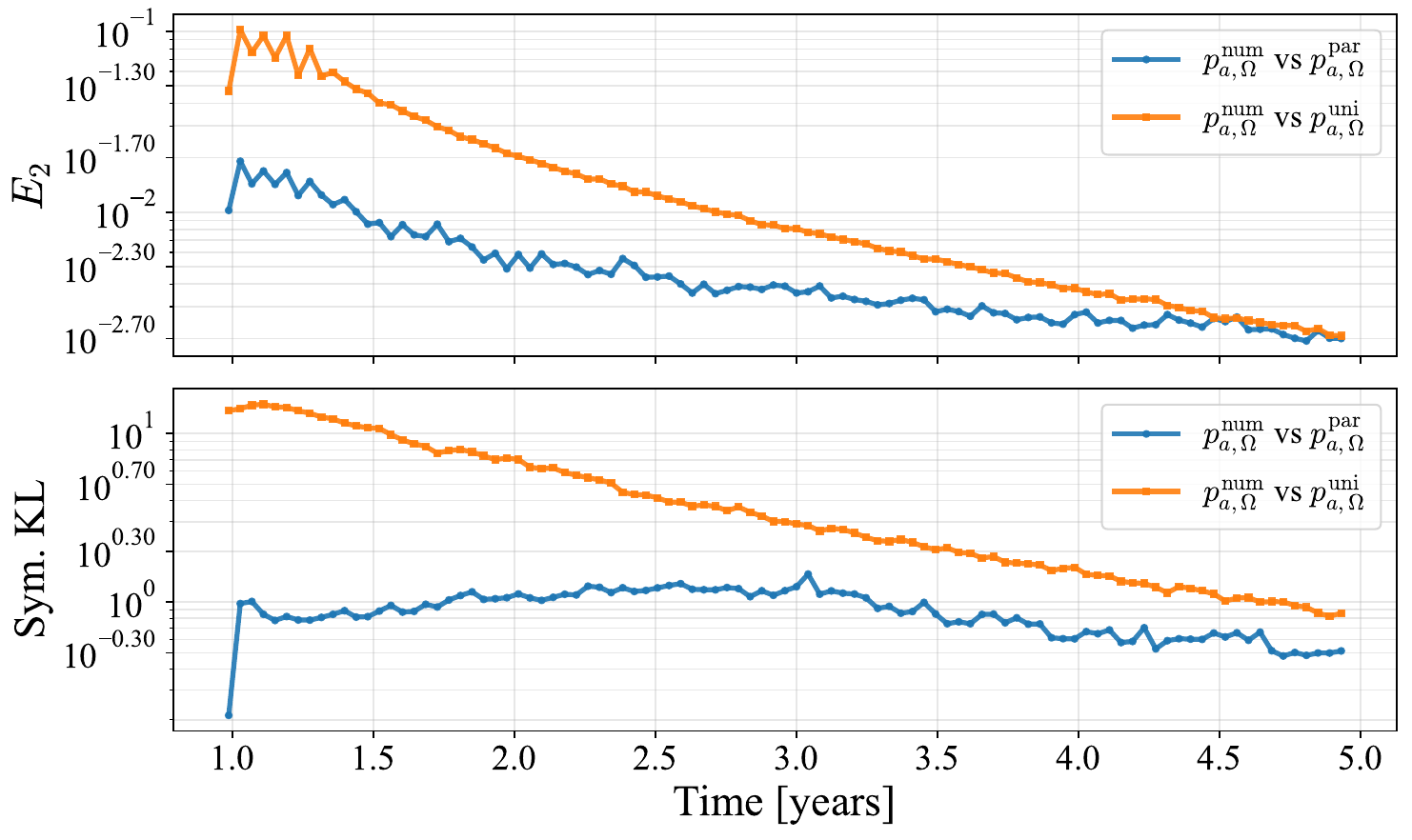}
\caption{$(a,\Omega)$}
\label{fig:e2kl_a_raan_log}
\end{subfigure}
\hfill
\begin{subfigure}[t]{0.48\textwidth}
\centering
\includegraphics[width=\linewidth]
{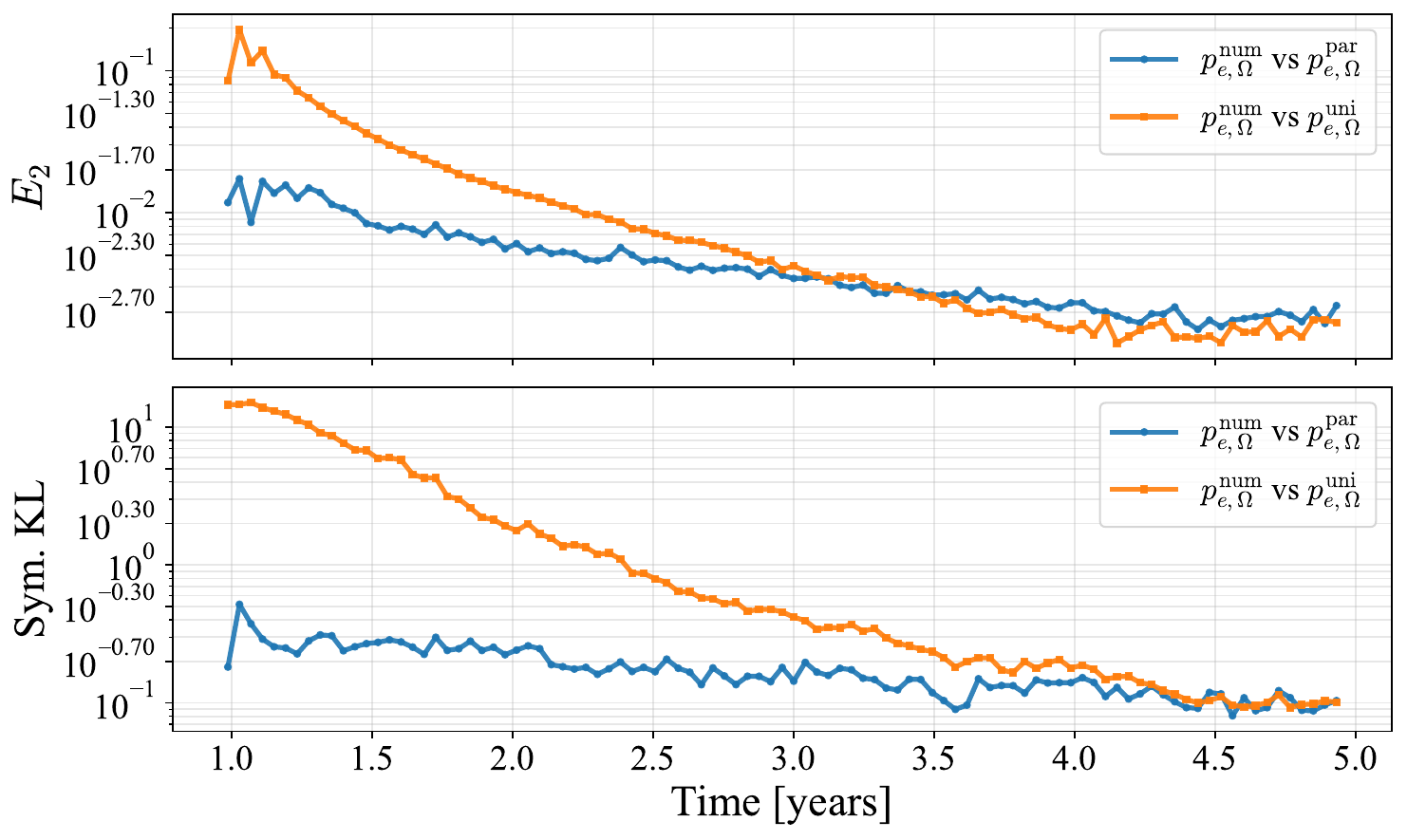}
\caption{$(e,\Omega)$}
\label{fig:e2kl_e_raan_log}
\end{subfigure}

\vspace{0.5em}

\begin{subfigure}[t]{0.60\textwidth}
\centering
\includegraphics[width=\linewidth]
{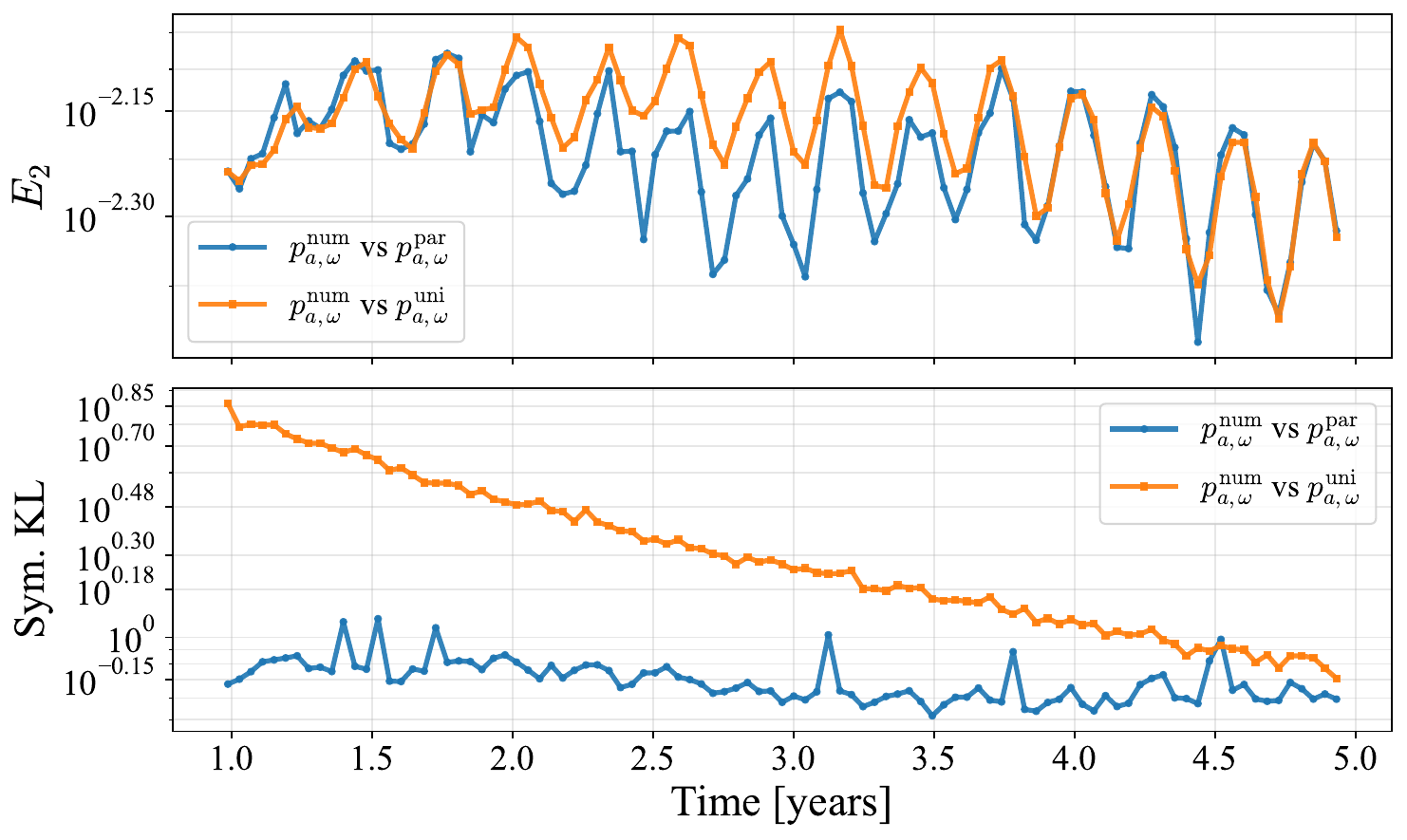}
\caption{$(a,\omega)$}
\label{fig:e2kl_e_aop_log}
\end{subfigure}
\caption{Distributional errors of the proposed surrogate and uniform-band approximation.}
\label{fig:e2kl_all}
\end{figure*}

The distributional accuracy of the proposed HGDM was assessed using the second-order error $E_2$ and the symmetric Kullback--Leibler divergence $D_{\mathrm{KL}}^{\mathrm{sym}}$ (see Appendix \ref{app:metrics}). The same metrics were computed for a baseline model corresponding to the classical uniform-band approximation.
Figure~\ref{fig:e2kl_all} shows that the proposed surrogate consistently outperforms the baseline for all considered distributions. The largest improvements are observed for pairs involving the \gls{raan}, namely $(a,\Omega)$ and $(e,\Omega)$, where the uniform-band approximation neglects the correlations induced by differential nodal precession.
Although the differences decrease as the cloud becomes more dispersed, the proposed surrogate remains systematically closer to the reference distribution throughout the entire propagation interval.

\subsection{Computational Efficiency and Storage Requirements}

Since the computational cost of the proposed framework is dominated by the propagation of the training fragments, a sensitivity analysis was performed by varying the training-set size $N_{\mathrm{train}}$. Table~\ref{tab:ntrain_kl} reports the average Kullback--Leibler divergence between the surrogate and reference distributions over a one-year propagation horizon. The error decreases rapidly with increasing $N_{\mathrm{train}}$ and reaches a plateau around $N_{\mathrm{train}}=1500$, indicating that accurate cloud reconstruction can be achieved from a relatively small number of propagated fragments.
\begin{table}[h]
\centering
\caption{\small Average KL divergence over the first year as a function of the training-set size.}
\label{tab:ntrain_kl}
\begin{tabular}{c|cccccccccc}
\hline
$N_{\mathrm{train}}$
& 250
& 500
& 1000
& 1500
& 2500
& 5000\\
\hline
$\overline{\mathrm{KL}}$
& 19.715
& 1.506
& 1.353
& 1.066
& 1.056
& \textbf{1.053}
\\
\hline
\end{tabular}
\end{table}
The corresponding runtimes are reported in Table~\ref{tab:runtime}. Relative to the reference Monte Carlo simulation, the proposed surrogate achieves speed-up factors of approximately $160\times$ and $54\times$ for $N_{\mathrm{train}}=500$ and $N_{\mathrm{train}}=1500$, respectively\footnote{Fragment propagations were performed on the CNES high-performance computing platform, whereas parameter estimation and cloud generation were executed on a standard desktop computer equipped with an Apple M2 processor (8 CPU cores) and 8 GB of RAM.}. These runtimes include fragment propagation, parameter estimation, and cloud generation.

\begin{table}[h]
\centering
\caption{Computation time for the LeoFrag scenario.}
\label{tab:runtime}
\begin{tabular}{lcc}
\hline
Method & Number of propagated fragments & Runtime \\
\hline
Reference Monte Carlo & $10^5$ & 9.26 h \\
Proposed model & $N_{\mathrm{train}}=1500$ & 10.36 min \\
Proposed model & $N_{\mathrm{train}}=500$ & 3.47 min \\
\hline
\end{tabular}
\end{table}
The reduction in storage requirements is also significant. Storing the fitted HGDM parameters requires approximately $1.2$ MB, compared with $1.4$ GB for the propagated debris-cloud ephemerides, corresponding to a compression factor exceeding three orders of magnitude \cite{Paul2025ICAWFragmentation}.

Overall, the proposed framework reconstructs high-fidelity fragmentation-cloud distributions from only a few hundred to a few thousand propagated fragments, while reducing computational cost by more than two orders of magnitude and storage requirements by more than three orders of magnitude.


\section{Conclusion}
This paper introduced a Hierarchical Generative Density Model (HGDM) for the long-term propagation of orbital fragmentation clouds. By exploiting the semi-major axis as the principal driver of cloud evolution and modeling the remaining orbital elements through conditional distributions, the proposed framework provides a compact probabilistic representation of the propagated debris-cloud distribution.
The method was evaluated on three fragmentation scenarios and compared against high-fidelity Monte Carlo reference simulations. The results showed that the proposed surrogate accurately reproduces the dominant multidimensional structures of fragmentation clouds, including persistent dependencies between semi-major axis, eccentricity, RAAN, and argument of perigee. In particular, it consistently outperformed classical band-formation approximations based on independent angular variables.
Beyond its accuracy, the proposed framework substantially reduces computational and storage requirements. Accurate cloud reconstructions were obtained using only a limited number of propagated representative fragments, yielding runtime reductions exceeding two orders of magnitude and storage reductions exceeding three orders of magnitude relative to conventional Monte Carlo propagation.
Future work will focus on continuous-time parameterizations of the model and on coupling the proposed density representation with large-scale debris-environment and collision-cascade simulations. In particular, the framework may provide an efficient foundation for studying fragmentation cascades and long-term debris growth associated with the Kessler syndrome.

\section*{Acknowledgments}
The author gratefully acknowledges the financial support from the École Normale Supérieure de Lyon (ENS Lyon) through the doctoral research contract \textit{N° ISAE: 2022-CIF-R-90}.

The first author gratefully acknowledges the French Space Agency (CNES) for hosting him during his research activities and for providing access to a unique operational environment in the field of Space Surveillance. Special thanks are extended to Valentin Barral, Head of the Space Surveillance Office at CNES for making this collaboration possible.

\bibliography{bibliography.bib}

\appendix

\section{Ariane~2~R/B Fragmentation Scenario}
\label{app:ariane}
\begin{figure}[H]
    \centering
    \begin{subfigure}[b]{0.49\textwidth}
        \centering
        \includegraphics[width=\textwidth]{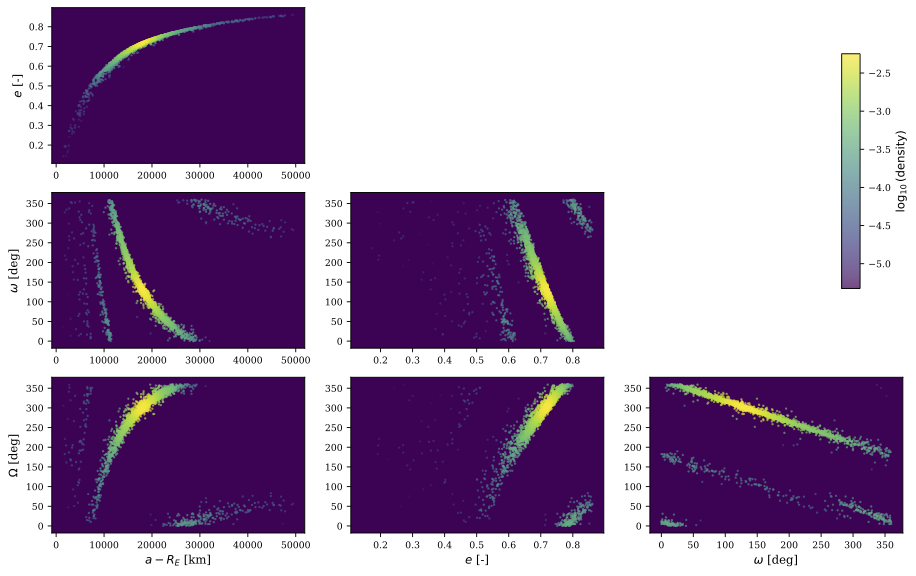}
        \caption{1 year after fragmentation.}
    \end{subfigure}
    \begin{subfigure}[b]{0.49\textwidth}
        \centering
        \includegraphics[width=\textwidth]{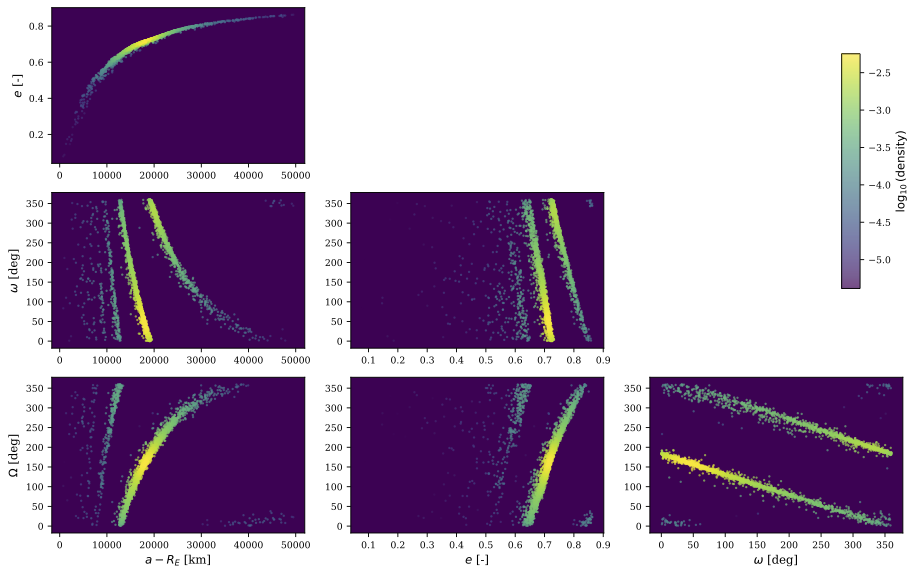}
        \caption{2 years after fragmentation.}
    \end{subfigure}
    \\  
    \begin{subfigure}[b]{0.49\textwidth}
        \centering
        \includegraphics[width=\textwidth]{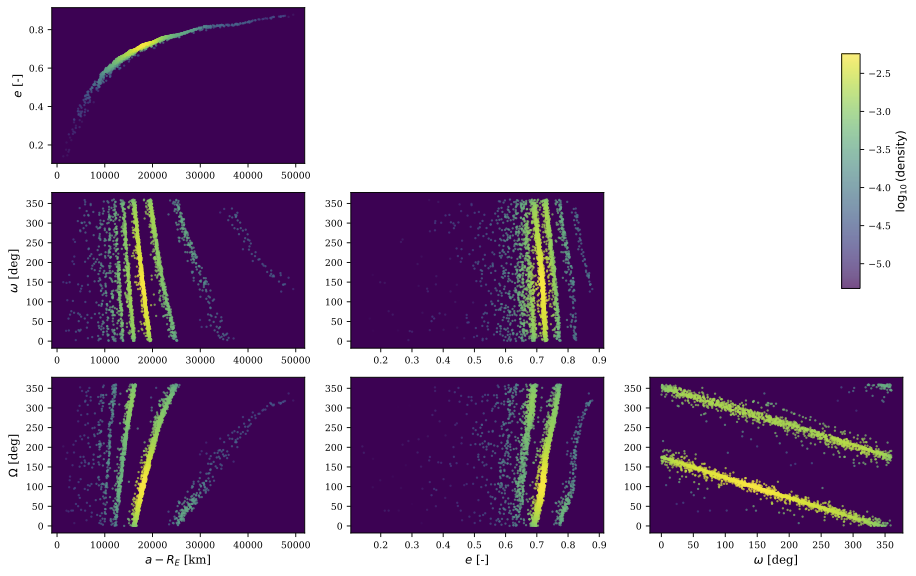}
        \caption{5 years after fragmentation.}
    \end{subfigure}
    \begin{subfigure}[b]{0.49\textwidth}
        \centering
        \includegraphics[width=\textwidth]{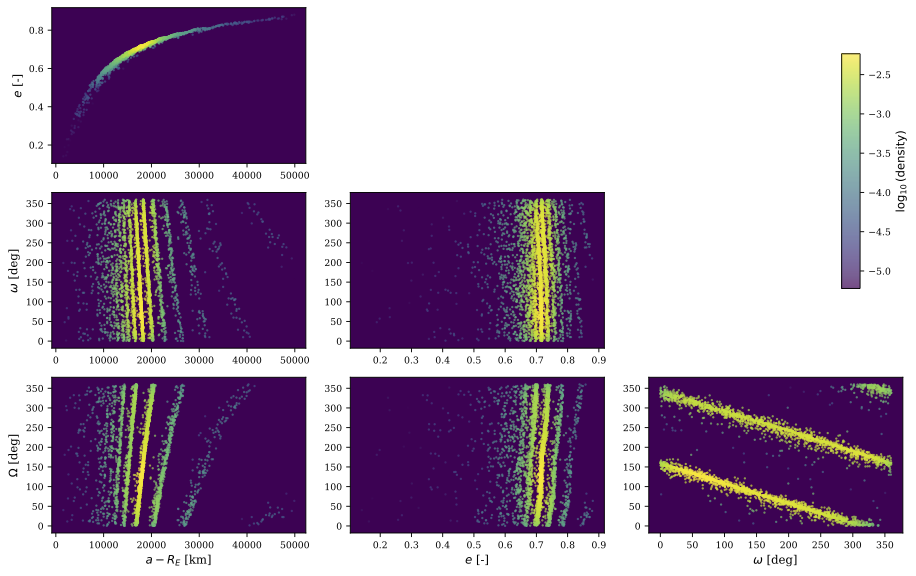}
        \caption{10 years after fragmentation.}
    \end{subfigure}
    \caption{\textbf{Ariane 2} scenario --- numerical reference
     at 1, 2, 5 and 10 ~years after fragmentation.}
    \label{fig:ariane2_ref}
\end{figure}
This appendix presents additional results obtained for the Ariane~2~R/B fragmentation scenario. Compared with the LeoFrag and ENVISAT cases discussed in the main text, this scenario corresponds to a highly eccentric orbit and exhibits a significantly broader dispersion in semi-major axis. As a consequence, the long-term evolution of the debris cloud is characterized by more complex dynamical structures and stronger nonlinear effects.

\begin{figure}[H]
    \centering
    \begin{subfigure}[b]{0.49\textwidth}
        \centering
         \includegraphics[width=\textwidth]{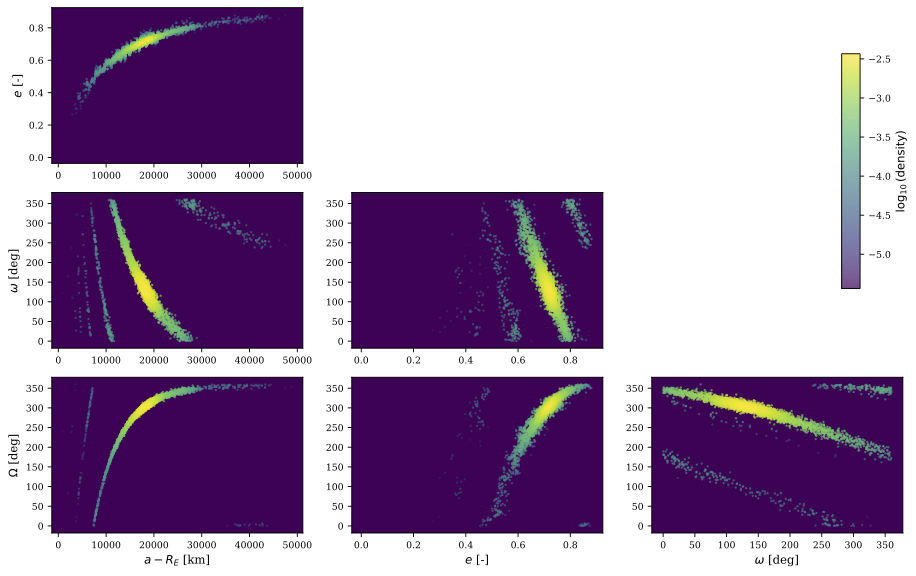}
        \caption{1 year after fragmentation.}
    \end{subfigure}
    \begin{subfigure}[b]{0.49\textwidth}
        \centering
        \includegraphics[width=\textwidth]{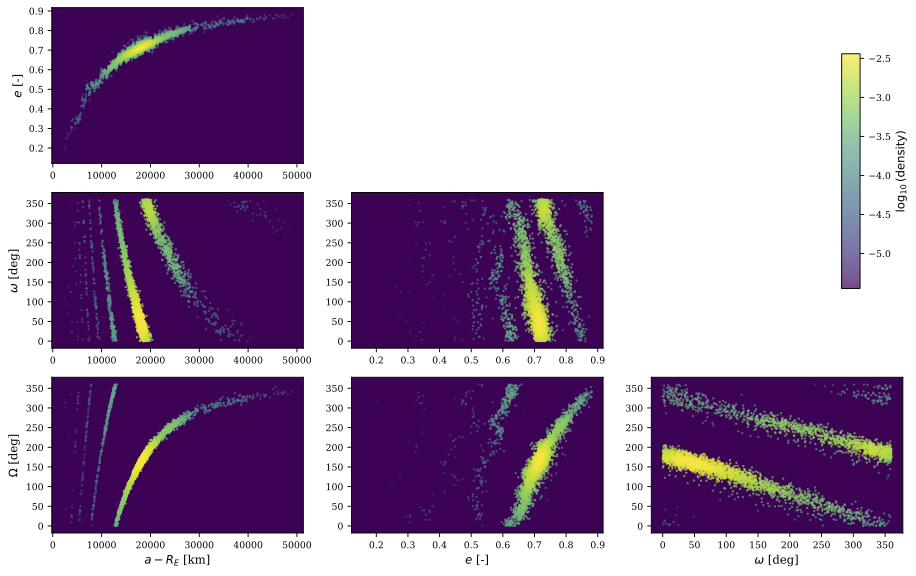}
        \caption{2 years after fragmentation.}
    \end{subfigure}
    \\  
    \begin{subfigure}[b]{0.49\textwidth}
        \centering
        \includegraphics[width=\textwidth]{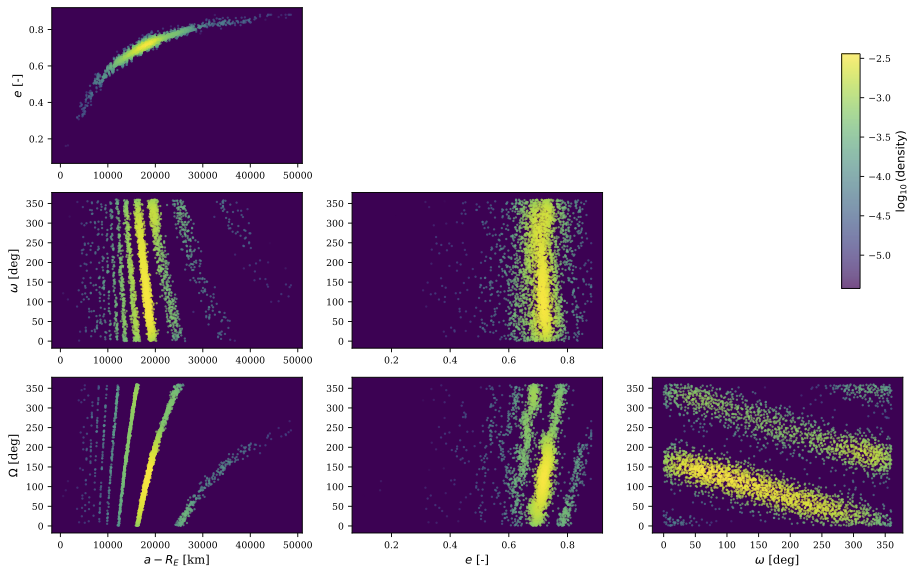}
        \caption{5 years after fragmentation.}
    \end{subfigure}
    \begin{subfigure}[b]{0.49\textwidth}
        \centering
        \includegraphics[width=\textwidth]{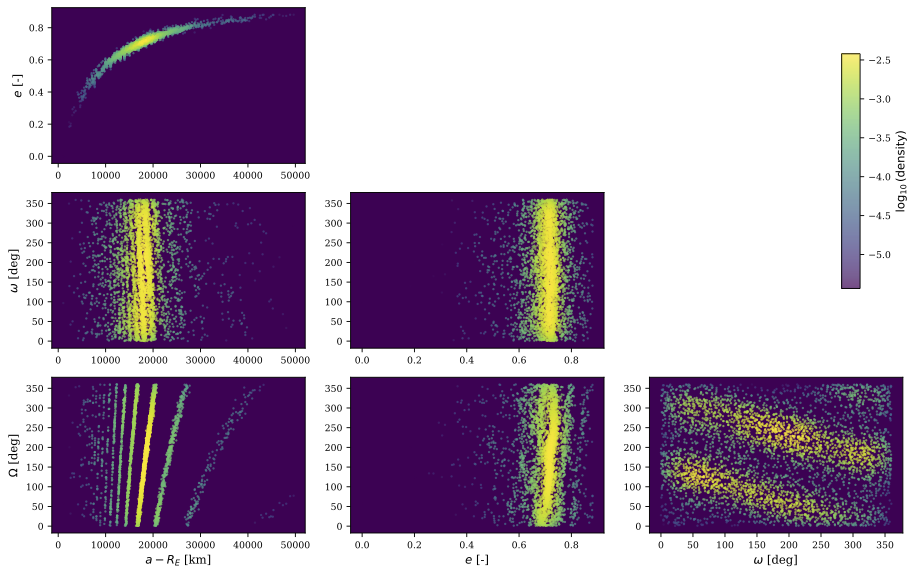}
        \caption{10 years after fragmentation.}
    \end{subfigure}
    
    \caption{\textbf{Ariane 2} scenario --- proposed surrogate model
     at 1, 2, 5 and 10 ~years after fragmentation.}    
     \label{fig:ariane2_pp}

\end{figure}

Figures~\ref{fig:ariane2_ref} and \ref{fig:ariane2_pp} compare the numerical reference distributions obtained from high-fidelity Monte Carlo propagation with the corresponding distributions generated by the proposed HGDM surrogate. Pairwise projections of the orbital-element space are shown after 1, 2, 5, and 10 years of propagation.
During the first years following fragmentation, the surrogate reproduces the principal features of the reference cloud with good accuracy. The dominant correlations involving the semi-major axis are preserved, and the overall morphology of the cloud remains consistent with the Monte Carlo reference. As the propagation duration increases, the cloud progressively develops more intricate structures, making accurate density reconstruction increasingly challenging. Consequently, discrepancies between the surrogate and the reference distributions become more visible at later epochs.

Despite these limitations, the HGDM continues to capture the dominant statistical trends governing the cloud evolution. These results indicate that the proposed hierarchical representation remains applicable beyond near-circular low-Earth-orbit fragmentation scenarios and can provide useful approximations for more complex orbital regimes while maintaining the substantial computational and storage advantages discussed in Section~\ref{sec:results}.

\section{Distributional Error Metrics}
\label{app:metrics}

To quantify the discrepancy between the reference distribution
$p^{\mathrm{ref}}$ and an approximate distribution $p$, two complementary metrics are employed.
The first metric is the second-order error
\begin{equation}
E_2(p,p^{\mathrm{ref}})=
\left(
\int
\left[
p(x)-p^{\mathrm{ref}}(x)
\right]^2
\,dx
\right)^{1/2},
\end{equation}
which measures the overall difference between the two densities in an $L^2$ sense.
The second metric is the symmetric Kullback--Leibler divergence,
\begin{equation}
D_{\mathrm{KL}}^{\mathrm{sym}}
(p,p^{\mathrm{ref}})=
D_{\mathrm{KL}}
(p\,|\,p^{\mathrm{ref}})
+
D_{\mathrm{KL}}
(p^{\mathrm{ref}}\,|\,p),
\end{equation}
where
\begin{equation}
D_{\mathrm{KL}}
(p\,|\,q)=
\int
p(x)
\log
\left(
\frac{p(x)}{q(x)}
\right)
dx .
\end{equation}
The metric $E_2$ quantifies absolute differences between densities, whereas $D_{\mathrm{KL}}^{\mathrm{sym}}$ measures discrepancies from an information-theoretic perspective. Lower values indicate a closer agreement with the reference distribution.

\end{document}